\documentclass[11pt]{amsart}
\usepackage[utf8]{inputenc}
\usepackage[english]{babel}

\usepackage{amscd,amssymb}
\usepackage{amsmath}
\usepackage{amsthm}
\usepackage{graphicx}
\usepackage{subfigure,fancybox}
\usepackage{epic,eepic}
\usepackage{amstext}
\usepackage{esint}
\usepackage[square,numbers]{natbib}
\usepackage{booktabs}
\usepackage[hide links]{hyperref}
\usepackage[noabbrev, capitalise, nameinlink]{cleveref}
\crefname{equation}{}{}
\crefname{assumption}{Assumption}{Assumptions}
\crefformat{equation}{\textup{#2(#1)#3}}
\usepackage{comment}
\usepackage{mathtools}
\usepackage{tikz,pgfplots}
\usepackage{stmaryrd}
\usetikzlibrary{plotmarks}
\usetikzlibrary{arrows.meta}
\usepgfplotslibrary{patchplots}
\usepackage{grffile}
\usepackage{tabularx}
\usepackage{float}

\def\lamTj{\Theta_{T,j}}
\def\bTj{b_{T,j}}

\newcommand{\Ca}{C_\mathrm{a}} 
\newcommand{\Cd}{C_\mathrm{d}} 
\newcommand{\Cbo}{C_\mathrm{bo}} 
\newcommand{\Cb}{C_\mathrm{b}} 
\newcommand{\Cla}{C_\Lambda} 
\newcommand{\Ct}{C_\mathrm{t}} 
\newcommand{\Col}{C_\mathrm{ol}} 
\newcommand{\Cp}{C_\mathrm{p}} 
\newcommand{\Ceo}{C}
\newcommand{\Cet}{C}
\newcommand{\Cdo}{C}
\newcommand{\Cdt}{C}
\newcommand{\Cslod}{C_\mathrm{slod}}
\newcommand{\Cloc}{C_\mathrm{lo}}
\newcommand{\Clod}{C_\mathrm{lod}}

\def\elem{T}

\def\V{\mathcal{V}}

\def\W{\mathcal{W}}

\def\P{\mathcal{P}}
\def\Pp{\mathcal{P}}

\newcommand{\N}{\mathbb{N}}

\def\ellp{m}

\def\T{\mathcal{T}}
\def\TH{\mathcal{T}_H}
\def\Th{\mathcal{T}_h}
\newcommand{\NH}{\mathcal{N}_{H}}

\newcommand{\with}{\,\colon\,}

\def\C{\mathcal{C}}
\def\Cl{\mathcal{C}^\ell}

\def\Clzt{\tilde{\mathcal{C}}^\ell_z}
\def\Ct{\mathcal{C}_T}
\def\Clt{\mathcal{C}_T^\ell}
\def\PiH{\Pi_H}

\def\BH{\mathcal B_H}
\def\diam{\mathrm{diam}}
\def\linh{\mathrm{span}}
\def\supp{\operatorname{supp}}

\newcommand{\Aop}{\mathcal{A}}

\newcommand{\abs}[1]{\lvert#1\vert}

\newcommand{\tnorm}[2]{{\left\lVert #1 \right\rVert}_{L^2(#2)}}
\newcommand{\tnormf}[2]{\| #1 \|_{L^2(#2)}}

\newcommand{\vnorm}[2]{{\left\lVert #1 \right\rVert}_{a,#2}}
\newcommand{\vnormf}[2]{\| #1 \|_{a,#2}}
\newcommand{\vnormof}[1]{\| #1 \|_{a,\Omega}}

\newcommand{\tspf}[3]{{( #1\,,\,#2 )}_{L^2(#3)}}

\newcommand\dx{\,\text{d}x}

\newcommand{\patch}{\mathsf{N}}
\newcommand{\patchTm}{\omega^\ellp_T}
\newcommand{\patchTmp}{\tilde{\omega}^\ellp_T}

\newcommand{\Vpumzi}{\V_{H,z}^{\ell}}

\newcommand{\Vpumz}{\V_{H,z}^{\ell,n}}
\newcommand{\Vpum}{\V_{H}^{\ell,n}}
\newcommand{\vpum}{u_{H}^{\ell,n}}

\newcommand{\Vslod}{\V_{H}^{\ell,\,\mathrm{SLOD}}}
\newcommand{\Vlod}{\V_{H}^{\ell,\,\mathrm{LOD}}}
\newcommand{\vslod}{u_{H}^{\ell,\,\mathrm{SLOD}}}
\newcommand{\vlod}{u_{H}^{\ell,\,\mathrm{LOD}}}
\newcommand{\vslodp}{u_{H}^{\ellp,\,\mathrm{SLOD}}}
\newcommand{\vlodp}{u_{H}^{\ellp,\,\mathrm{LOD}}}

\newcommand{\vslodpz}{u_{H,z}^{\ellp,\,\mathrm{SLOD}}}
\newcommand{\vlodpz}{u_{H,z}^{\ellp,\,\mathrm{LOD}}}
\newcommand{\pl}{\psi_{T,j}^{\ell}}

\newcommand{\g}{g_{T,j}^{\ell}}
\newcommand{\plp}{\psi_{T,j}^{\ellp}}
\newcommand{\ppt}{\tilde\varphi_{T,j}^{\ellp}}
\newcommand{\plpt}{\tilde\psi_{T,j}^{\ellp}}
\newcommand{\plph}{\widehat{\psi}_{T,j}^{\ellp}}
\newcommand{\gp}{g_{T,j}^{\ellp}}
\newcommand{\gpt}{\tilde g_{T,j}^\ellp}

\definecolor{myBlue}{RGB}{113,104,238} 
\definecolor{myGreen}{RGB}{114,175,30} 
\definecolor{myRed}{RGB}{180,50,50}  
\definecolor{myOrange}{RGB}{225,92,22} 

\definecolor{color0}{rgb}{0.65,0,0.15}
\definecolor{color1}{rgb}{0.84,0.19,0.15}
\definecolor{color2}{rgb}{0.96,0.43,0.26}
\definecolor{color3}{rgb}{0.99,0.68,0.38}
\definecolor{color4}{rgb}{1,0.65,0.56}
\definecolor{color5}{rgb}{0.67,0.85,0.91}
\definecolor{color6}{rgb}{0.27,0.46,0.71}

\def\SLODColor{color0}
\def\PUMnOne{color4}
\def\PUMnTwo{color3}
\def\PUMnThree{color1}
\def\PUMnFour{color6}
\def\PUMnFive{color5}

\def\lOneMarker{square*}
\def\lTwoMarker{triangle*}
\def\lThreeMarker{diamond*}
\def\pMarker{*}

\def\dx{\,\text{d}x}

\newcommand{\dnz}{d_n^z}
\newcommand{\dn}{d_n}

\newtheorem{theorem}{Theorem}[section]

\newtheorem{assumption}[theorem]{Assumption}

\theoremstyle{definition}

\theoremstyle{remark}
\newtheorem{remark}[theorem]{Remark}
\numberwithin{theorem}{section}
\numberwithin{equation}{section}
\numberwithin{table}{section}
\numberwithin{figure}{section}

\allowdisplaybreaks

\usepackage[colorinlistoftodos,prependcaption]{todonotes}
\begin{document}
	
\begin{abstract}
	This paper presents a novel multi-scale method for elliptic partial differential equations with arbitrarily rough coefficients. In the spirit of numerical homogenization, the method constructs problem-adapted ansatz spaces with uniform algebraic approximation rates. Localized basis functions with the same super-exponential localization properties as the recently proposed Super-Localized Orthogonal Decomposition enable an efficient implementation. The method's basis stability is enforced using a partition of unity approach. A natural extension to higher order is presented, resulting in higher approximation rates and enhanced localization properties. 
	We perform a rigorous a priori and a posteriori error analysis and confirm our theoretical findings in a series of numerical experiments. In particular, we demonstrate the method's applicability for challenging high-contrast channeled coefficients. 
\end{abstract}
	
\title[A super-localized generalized finite element method]{A super-localized generalized finite element method}

\author[P.~Freese, M.~Hauck, T.~Keil, D.~Peterseim]{Philip Freese$^\dagger$, Moritz Hauck$^\dagger$, Tim Keil$^\ddagger$, Daniel Peterseim$^*$}
\address{${}^{\dagger}$ Institute of Mathematics, University of Augsburg, Universit\"atsstr.~12a, 86159 Augsburg, Germany}
\email{\{philip.freese, moritz.hauck\}@uni-a.de}
\address{${}^{*}$ Institute of Mathematics \& Centre for Advanced Analytics and Predictive Sciences (CAAPS), University of Augsburg, Universit\"atsstr.~12a, 86159 Augsburg, Germany}
\email{daniel.peterseim@uni-a.de}
\address{${}^{\ddagger}$ Mathematics Münster, Westfälische Wilhelms-Universität Münster, Einsteinstr. 62, 48149 M\"unster, Germany}
\email{tim.keil@wwu.de}
\thanks{The work of Philip Freese, Moritz Hauck, and Daniel Peterseim is part of a project that has received funding from the European Research Council ERC under the European Union's Horizon 2020 research and innovation program (Grant agreement No.~865751).}
\thanks{Tim Keil acknowledges funding by the Deutsche Forschungsgemeinschaft under Germany’s Excellence Strategy EXC 2044 390685587, Mathematics M\"unster: Dynamics -- Geometry -- Structure.}
\maketitle

\indent
{\tiny {\bf AMS subject classifications.}  
	{\bf 65N12}, 
	{\bf 65N30} 
} 

\section{Introduction}

We consider the numerical solution of a second-order linear elliptic partial differential equation with a strongly heterogeneous coefficient. The coefficient may be non-periodic, with oscillations appearing on several non-separated scales. For such coefficients, classical finite element methods based on problem-independent polynomial ansatz spaces typically yield unsatisfactory approximations, cf.~\cite{BabO00}. It is possible to overcome this issue by incorporating problem-specific information into the method's ansatz space, which is commonly known under the term numerical homogenization and has been an active research field throughout the past decades. For an overview on numerical homogenization, we refer to the recent textbooks \cite{OwhS19,MalP20} and the review article~\cite{Peterseim2021}.

Under minimal assumptions on the coefficient, numerical homogenization is able to achieve optimal orders of approximation without any pre-asymptotic effects. However, this is not possible without a computational overhead. Compared to classical finite element methods, it is either necessary to consider basis functions with an enlarged support, or to increase the number of basis functions per mesh entity. As prominent examples, we mention the Multiscale Spectral Generalized Finite Element Method~(MS-GFEM)~\cite{BaL11,EFENDIEV2013116,Ma22}, the Adaptive Local Basis~\mbox{(AL-Basis)}~\cite{GGS12,Wey17}, the Localized Orthogonal Decomposition method~(LOD)~\cite{HeP13,MaP14,KPY18,BrennerLOD}, Rough Polyharmonic Splines~(RPS)~\cite{OZB14}, and gamblets~\cite{Owh17}. 

The above-mentioned approaches can be distinguished into two classes. Methods like MS-GFEM and the AL-Basis first solve local spectral problems in the space of (locally) operator-harmonic functions. The respective ansatz spaces are then constructed by gluing together local eigenfunctions by means of a partition of unity~\cite{BaMe96,BaM97}.
For such methods, the support of the basis functions is fixed by the choice of partition of unity. For convergence of optimal order, the number of local eigenfunctions taken into account needs to be increased logarithmically with the desired accuracy. In order to make these approaches computationally more efficient, one can use random sampling strategies as proposed, e.g., in~\cite{BuS18}. 

The second class of methods includes the LOD, RPS, and gamblets. The idea is to construct problem-adapted ansatz spaces by applying the solution operator to specific classical finite element spaces with respect to some coarse mesh~$\TH$, which typically does not resolve the coefficient's oscillations. Due to their connection to isogeometric analysis in the case of constant coefficients, such methods are sometimes referred to as spline-type approaches.
While optimal approximation orders of these methods are achieved by design, the true challenge is to construct a local basis of the problem-dependent ansatz space. 
An almost optimal solution is provided by the LOD, which constructs a fixed number of basis functions per mesh entity that decay exponentially fast with respect to the coarse mesh. This rapid decay enables a localization of the basis functions to \mbox{$\ell$-th} order element patches with diameters of order~$\ell H$. For convergence of optimal order, the oversampling parameter~$\ell$ needs to be increased logarithmically with the desired accuracy. 

Recently, the Super-Localized Orthogonal Decomposition (SLOD) has been proposed in~\cite{HaPe21b} (see also \cite{Freese-Hauck-Peterseim,Bonizzoni-Freese-Peterseim,graphSLOD}). The key contribution of the SLOD is a novel localization strategy that enables a significantly improved localization compared to the LOD. The SLOD constructs rapidly decaying basis functions, yielding super-exponentially decaying localization errors. These improved localization leads to smaller local patch problems for the basis computation and a sparser coarse system matrix. Numerical experiments indicate that the SLOD outperforms the LOD, achieving similar magnitude errors for significantly smaller oversampling parameters. Until now, for the best practical realization of the SLOD in~\cite{HaPe21b}, the stability of the SLOD basis functions cannot be guaranteed a priori. 
For high-contrast channeled coefficients or convection-dominated regimes, these basis stability issues may deteriorate the method's approximation quality (see the numerical experiments in~\cref{sec:Implementation and numerical experiments}).

This paper proposes a novel multi-scale method that, on the one hand, preserves the unique localization properties of the SLOD and, on the other hand, resolves the aforementioned basis stability issues. This is achieved by combining the SLOD with a partition of unity approach. More precisely, locally on nodal patches, we apply the respective local solution operator to classical finite element source terms.  Multiplying these spaces with the corresponding hat-functions yields local ansatz spaces with a low effective dimension. Consequently, low-dimensional optimally approximating spaces are constructed by solving local spectral problems. 
Compared to MS-GFEM methods, the proposed method has the major advantage that the local spectral problems are posed in a space spanned by a small number of deterministic snapshots. Hence, possible random sampling strategies can be avoided. Furthermore, due to their low dimension, the local spectral problems are easy to solve. The  global problem-adapted ansatz space is then obtained by gluing together the local optimal approximation spaces using a partition of unity. We highlight that, from an application point of view, the proposed multi-scale method is conceptually simple and straightforward to implement. Further, by adapting the polynomial degree of the finite element source terms, one can easily construct higher-order versions of the method. Similarly as for the higher-order LOD \cite{Mai20ppt,DHM22}, one obtains higher-order convergence rates using the regularity of the source only.

We prove that the proposed method possesses the advantageous localization and convergence properties of the SLOD which can be quantified a posteriori. Building on the well-understood theoretical foundation of the LOD, we additionally perform a pessimistic a priori error analysis, proving that the proposed method at least recovers the convergence and localization properties of the LOD. For the method's higher-order versions, we observe that solely increasing the polynomial degree significantly improves the localization properties. Numerical experiments even suggest that an (almost) local basis exists for sufficiently large polynomial degrees.
Another noteworthy contribution of this work is the implementation of the proposed method as well as the SLOD in the \texttt{Python}-library \texttt{gridlod} \cite{gridlod}. This serves the principles of open access and reproducibility while enabling computations on large parallel clusters.

The outline of this paper is as follows. In \cref{sec:Model problem}, we introduce the prototypical elliptic model problem. \Cref{sec:Preliminaries} recalls preliminary results, which are then used in \cref{sec:Multi-scale method} for the definition of the novel multi-scale method. An a posteriori error analysis is presented in \cref{sec:a posteriori error analysis} followed by a pessimistic a priori error analysis in \cref{sec:a priori error analysis}. Finally,  \cref{sec:Implementation and numerical experiments} presents a series of numerical experiments which confirm our theoretical findings.

\section{Model problem} \label{sec:Model problem}

We consider the prototypical second-order elliptic PDE $-\mathrm{div}A\nabla u = f$ in weak form with homogeneous Dirichlet boundary conditions on a polygonal/polyhedral Lipschitz domain $\Omega\subset \mathbb{R}^d$, $d \in\{2,3\}$.
Furthermore, without loss of generality, we assume that $\Omega$ is scaled such that its diameter is of order one.
The coefficient function $A\in L^\infty(\Omega,\mathbb{R}^{d\times d})$ may be matrix-valued and is assumed to be symmetric and positive definite almost everywhere. More specifically, we assume that there exist constants $0<\alpha\leq \beta< \infty$ such that
\begin{equation}\label{eq:propA}
	\alpha \abs{\eta}^2 \leq (A(x)\eta)\cdot \eta \leq \beta \abs{\eta}^2,\qquad x \in \Omega,\;\eta\in \mathbb{R}^d
\end{equation}
with $\abs{\cdot}$ denoting the Euclidean norm. The weak formulation of the elliptic model problem uses the Sobolev space $\V \coloneqq H^1_0(\Omega)$ and the bilinear form  $a\colon\V\times\V\rightarrow \mathbb{R}$, given by
\begin{equation*}
	a(u, v) \coloneqq
	\int_\Omega (A\nabla u)\cdot\nabla v\dx.
\end{equation*}
The symmetry and the condition \cref{eq:propA} ensure that the above bilinear form is an inner product on $\V$. Its induced norm is the energy norm $\vnorm{\cdot}{\Omega} \coloneqq \sqrt{a(\cdot , \cdot)}$, which is equivalent to the canonical Sobolev norm on $\V$. The Lax--Milgram theorem ensures that, for all source terms $f\in L^2(\Omega)$, there exists a unique weak solution $u \in \V$ to the boundary value problem, satisfying
\begin{equation}\label{eq:sol}
	a(u, v) = (f,v)_{L^2(\Omega)},\qquad\text{for all }v \in \V.
\end{equation}
Note that the moderate restriction to source terms in $L^2(\Omega)$ (rather than the dual space  $\V^\prime=H^{-1}(\Omega)$) will be essential for the uniform convergence of the numerical homogenization method. We note that the possibly rough coefficient generally prevents $H^2(\Omega)$-regularity of the solution, which would be required by classical finite elements. For 
 a generalization to source terms with less regularity, we refer to~\cite{Peterseim2021}. 
Let us mention that the proposed method is not restricted to the class of elliptic PDEs with Dirichlet boundary conditions. Considering the extensions of the SLOD \cite{Freese-Hauck-Peterseim,Bonizzoni-Freese-Peterseim}, also an extension of the proposed method to (non-symmetric) coercive operators and even Helmholtz-type problems seems possible. In particular, more general boundary conditions of Neumann and Robin type can be taken into account. 

Henceforth, we refer to $\mathcal A^{-1}\colon L^2(\Omega)\rightarrow \V$ as the solution operator that maps $f\in L^2(\Omega)$ to the unique solution $u\in\V$ of \cref{eq:sol}.
Moreover, for a subdomain $\omega\subset \Omega$, we denote by $a_\omega(\cdot,\cdot)$ and $\vnormf{\cdot}{\omega}$ the restriction of the bilinear form $a(\cdot,\cdot)$ to $\omega$ and the restricted energy norm, respectively.
The restricted solution operator subject to homogeneous Dirichlet boundary conditions on $\partial \omega$ is denoted by $\Aop^{-1}_{\omega}\colon L^2(\omega)\rightarrow  H^1_0(\omega)$.

\section{Preliminaries} \label{sec:Preliminaries}
Let $\TH$ denote a quasi-uniform coarse mesh of $\Omega$ consisting of closed, simplicial, or quadrilateral shape-regular elements. 
The subscript~$H$ denotes the maximal element diameter, i.e., $H \coloneqq \max_{T\in\TH}\diam(T)$ and by~$\NH$, we denote the set of all (interior and boundary) vertices of~$\TH$. For the ease of presentation, we henceforth only consider quadrilateral meshes; the extension to triangular meshes is straightforward.

The proposed multi-scale method utilizes the concept of patches. For any $\ell\in\N_0$, we define the $\ell$-th order (element) patch of a union of elements $\omega\subset\Omega$ recursively by 
\begin{align*}
	\patch^0(\omega) \coloneqq \omega,\qquad \patch^{\ell + 1}(\omega) \coloneqq \bigcup\left\lbrace T\in\TH \with T \cap \patch^\ell(\omega) \neq \emptyset \right\rbrace.
\end{align*}

\subsection{Discontinuous finite element spaces}
For a fixed (but arbitrary) polynomial degree $p$, we denote with
\begin{equation*}
	\P(\TH) \coloneqq  \{v \in L^2(\Omega)\with v\vert_T, T \in \TH, \text{ is a polynomial of coordinate degree }\leq p\}
\end{equation*}
the non-conforming space (with respect to $\V$) consisting of element-wise defined polynomials. 
We define the restriction of $\P(\TH)$ to a subdomain $\omega \subset \Omega$ by
\begin{equation*}
	\P(\omega) \coloneqq  \{v \in \P(\TH)\with \supp(v) \subset \omega\}.
\end{equation*}
One can characterize the space $\P(T)$, $T \in \TH$ in terms of a suitable orthonormal basis $\{\lamTj\with j=1,\dots,J\}$ with $J\coloneqq \dim(\P(T))$, e.g., shifted tensor-product Legendre polynomials. Hence, a local orthonormal basis of $\P(\TH)$ is given by $\{\lamTj\with T \in \TH, j = 1,\dots,J\}$.

Let $\PiH\colon L^2(\Omega) \to \P(\TH)$ denote the $L^2$-orthogonal projection, which for each $v \in L^2(\Omega)$, is given by the element-wise equation
\begin{equation*}\label{eq:L2proj}
	(\PiH v, w)_{L^2(T)} = (v,w)_{L^2(T)}, \qquad\text{for all }w \in \P(T),\;T \in \TH.
\end{equation*}
The projection satisfies the following stability and approximation estimates
\begin{alignat}{2}
	\tnormf{\PiH v}{T}&\leq \tnormf{v}{T},\qquad &&v \in L^2(T), T\in \TH,\label{e:stabL2}\\
	\tnormf{(1-\PiH)v}{T}&\leq \Ca H \tnormf{\nabla v}{T},\qquad &&v \in H^1(T), T \in \TH
\end{alignat}
with constant $\Ca>0$ depending only on the regularity of the mesh $\TH$ and the polynomial degree $p$, see, e.g., \cite{HSS02}. 

In addition, we define the broken Sobolev space $H^k(\TH)$, $k \in \N$, by
\begin{equation*}
H^k(\TH) \coloneqq \big\{ v \in L^2(\Omega)\with v\vert_T \in H^k(T),\;T \in \TH \big\}
\end{equation*}
with the seminorm
\begin{equation*}
\abs{\cdot}^2_{H^k(\TH)} \coloneqq \sum_{T \in \TH} \abs{\cdot}^2_{H^k(T)}.
\end{equation*}
\subsection{Conforming companion spaces}
We next define local conforming companions of the functions~$\lamTj$, so-called bubble functions, which we denote by $\bTj$. For each element $T \in \TH$, these functions fulfill, for $1 \leq j \leq J$,
\begin{equation}\label{eq:bubble}
	\bTj \in H^1_0(T), \qquad \PiH \bTj = \lamTj. 
\end{equation}
We do not require an explicit characterization. However, it is important that such functions actually exist. This is guaranteed by~\cite[Cor.~3.6]{Mai20ppt} stating that, for all $\lamTj$, there exist a corresponding bubble function $b_{T,j}$ such that
\begin{equation}
	\label{eq:bubbleest}
	\tnormf{b_{T,j}}{T} + H\tnormf{\nabla b_{T,j}}{T} \leq \Cb \tnormf{\lamTj}{T},\qquad T\in \TH,j=1,\dots,J 
\end{equation}
with constant $\Cb>0$ depending solely on the mesh regularity of $\TH$ and the polynomial degree $p$. 

By means of the bubble functions, we can define the operator $\BH$ mapping possibly non-conforming functions to $\TH$-piecewise bubble functions with the same $L^2$-projection.
For any function in $\P(\TH)$, we uniquely define $\BH$ by setting $\BH \lamTj \coloneqq \bTj$ for all $T \in \TH$, $j =1,\dots,J$. We can extend the operator to $L^2(\Omega)$ by setting
\begin{equation*}
	\BH v \coloneqq \BH \PiH v,\qquad \text{for all }v \in L^2(\Omega).
\end{equation*} 
Clearly, the kernels of the operators $\PiH$ and $\BH$ coincide and one can prove the local stability estimate
\begin{equation}
	\tnormf{\BH v}{T} + H\tnormf{\nabla \BH v}{T} \leq \Cbo \tnormf{v}{T},\qquad v \in L^2(T), T\in\TH
	\label{eq:bubbleopest}
\end{equation}
with another constant $\Cbo>0$ depending solely on the mesh regularity of $\TH$ and the polynomial degree $p$. By the definition of $\BH$, we obtain, for all $v \in L^2(\Omega)$, $q \in \P(\TH)$,
\begin{equation}
	\label{eq:BHL2proj}
	\tspf{\BH v}{q}{\Omega} =  \tspf{\BH\PiH v}{q}{\Omega} = \tspf{\PiH v}{q}{\Omega} =  \tspf{q}{v}{\Omega},
\end{equation}
i.e., $\BH$ is the $L^2$-projection onto the space of bubble functions.

\subsection{Partition of unity}
The proposed multi-scale method is based on the framework of partition of unity methods, cf.~\cite{BaMe96,BaM97}. Although, there is great flexibility in the choice of such a partition, for simplicity, we restrict ourselves to the hat-functions $\{\Lambda_z\with z \in \NH\}$ corresponding to all (interior and boundary) nodes of $\TH$. Recall that the hat-function $\Lambda_z$ associated with node $z \in \NH$ is a continuous $\TH$-piecewise bilinear function uniquely defined by setting its nodal values for all $y \in \NH$ to $\Lambda_z(y) = \delta_{yz}$ with~$\delta$ denoting the Kronecker symbol. 
By definition, the hat-functions have an $L^\infty$-norm of one and, due to the shape-regularity of $\TH$, their gradients satisfy 
\begin{equation}
	\|\nabla \Lambda_z\|_{L^{\infty}(\Omega)} \leq \Cla H^{-1},\qquad z\in \NH\label{eq:boundhatfun}
\end{equation}
with constant $\Cla>0$ depending solely on the mesh regularity of $\TH$.
Denoting $\omega_z\coloneqq \operatorname{supp}(\Lambda_z)$, the shape-regularity of $\TH$ also implies that the supports $\{\omega_z\with z \in \NH\}$ have a finite overlap, i.e., the maximal number of overlapping supports 
\begin{equation}
	\Col \coloneqq \max_{T\in\TH}\#\{z \in \NH \with  T \subset \omega_z\}
	\label{eq:overlap}
\end{equation}
is uniformly bounded.
Subsequently, we abbreviate the node patches around $z \in \NH$ and the element patches around $T \in \TH$ by $\omega_z^\ell \coloneqq \patch^\ell(\omega_z)$ and  $\omega_T^\ell \coloneqq \patch^\ell(T)$, respectively.

\section{Multi-scale method} \label{sec:Multi-scale method}

This section introduces the proposed multi-scale method. The local ansatz spaces of the method are constructed by applying the local solution operator on an oversampling domain to piecewise polynomial source terms and by subsequent restriction to a subdomain.  
Due to the oversampling, the resulting local spaces have a low effective dimension, and thus, low-dimensional optimally approximating spaces are utilized.
The ansatz space of the method is obtained by gluing together the low-dimensional local approximation spaces. 
Note that we consider a fixed polynomial degree $p$ and do not track the dependence of constants on $p$. Explicitly tracking this dependence would make the analysis less clear and also add no value, as it relies on estimates that are pessimistic in $p$. 

\subsection{Local approximation spaces}
For any $z \in \NH$, we aim to approximate the restriction of the solution space $\V\lvert_{\omega_z}$ using local approximation spaces. This is accomplished by choosing the local approximation space $\Vpumzi\lvert_{\omega_z}$ with
 \begin{align}\label{eq:V_H_ell_z}
 \Vpumzi \coloneqq \linh\big\lbrace \Aop^{-1}_{\omega_z^\ell} \, q
  \with q \in \Pp(\omega_z^\ell) \big\rbrace \subset H^1_0(\omega_z^\ell)
 \end{align}
 being defined on the oversampling domain $\omega_z^\ell$.
 
  Due to the oversampling, the restricted space contains many redundant functions. This holds, in particular, after the multiplication with the hat-function $\Lambda_z$ when gluing the local approximation spaces together. Hence, we investigate the optimal approximation of $\Lambda_z\Vpumzi$ by $n$-dimensional subspaces $Q(n) \subset H^1_0(\omega_z)$. Given the subspace $Q(n)$, its worst-case best approximation error is defined as
\begin{equation*}
\adjustlimits \sup_{v \, \in \, \Vpumzi} \inf_{w \in Q(n)} \frac{\vnormf{\Lambda_zv -w}{\omega_z}}{\vnormf{v}{\omega_z^\ell}}.
\end{equation*}
Typically, the minimal worst-case best approximation error is referred to as Kolmogorov $n$-width, cf.~\cite{Pin85}, and is defined as 
\begin{equation}
\label{Kolmogorovnw}
\dnz(H,\ell)\coloneqq \adjustlimits \inf_{Q(n) \subset H^1_0(\omega_z)} \sup_{v \, \in \, \Vpumzi} \inf_{w \in Q(n)} \frac{\vnormf{\Lambda_zv -w}{\omega_z}}{\vnormf{v}{\omega_z^\ell}}.
\end{equation}
Indeed, there exists a corresponding optimal local approximation space of dimension~$n$, which we explicitly compute. For this, we solve the low-dimensional eigenvalue problem, which seeks eigenpairs $(v,\lambda) \in \Vpumzi\times \mathbb R$ such that
\begin{equation}\label{eq:evp}
a_{\omega_z}(\Lambda_zv, \Lambda_zw) = \lambda \, a_{\omega_z^\ell}(v,w),\qquad\text{for all }w \in \Vpumzi.
\end{equation}
We denote the eigenfunctions by $\{v_i\with i = 1,\dots,N\coloneqq \dim(\Vpumzi)\}$ assuming an ordering such that the corresponding eigenvalues satisfy $\lambda_1\geq \lambda_2\geq\dots\geq \lambda_N \geq 0$.
Consequently, denoting 
\begin{equation*}
\Vpumz \coloneqq \linh\{v_i\with i=1,\dots,n\},
\end{equation*}
the optimal local approximation space of dimension $n$ is given by $\Lambda_z\Vpumz$.

\subsection{Global approximation space}
A global approximation space is obtained by gluing together the above local approximation spaces using the partition of unity, i.e., 
\begin{align*}
\Vpum \coloneqq \sum_{z \in \NH}  \Lambda_z\Vpumz.
\end{align*}
We measure the overall error when approximating $\Lambda_z\Vpumzi$ by spaces of dimension $n$ by
\begin{equation}
\label{e:kolmogorov}
\dn(H,\ell) \coloneqq \max_{z\in\NH} \dnz(H,\ell).
\end{equation}
Finally, the proposed method seeks $\vpum \in \Vpum$ such that
\begin{equation} \label{eq:PUMSLOD}
	a(\vpum,v) = \tspf{f}{v}{\Omega},\qquad\text{for all }v \in \Vpum.
\end{equation}
Note that the mesh size $H$ determines the accuracy of the approximation. The oversampling parameter $\ell$ specifies the size of the local patch problems and determines the method's localization error, whereas the number of local functions is given by $n$ and needs to be chosen sufficiently large. For a precise choice of the parameters $H,\ell$, and $n$, we refer to \cref{rem:bound_delta_LOD,rem:bound_delta_SLOD}.
The following two sections are devoted to the theoretical analysis of the proposed multi-scale method. In \Cref{sec:a posteriori error analysis}, we present an a posteriori error analysis, while in \Cref{sec:a priori error analysis}, we present a priori error bounds.  

\section{A posteriori error analysis} \label{sec:a posteriori error analysis}
Subsequently, we derive an a posteriori error analysis of the proposed method by establishing a connection to the SLOD introduced in \cite{HaPe21b}. The SLOD is conceptually related to the proposed method, as it also constructs its basis functions by applying the local solution operator to $\TH$-piecewise polynomial source terms.

\subsection{Higher-order SLOD}
For the a posteriori error analysis, we first briefly introduce a higher-order variant of the SLOD. Note that this variant only serves theoretical purposes and is not investigated numerically.
Let us fix an arbitrary element $\elem\in\T_H$ and oversampling parameter $\ell\in\mathbb{N}$. Henceforth, we drop all fixed indices and denote the $\ell$-th order patch around $T$ just by $\omega \coloneqq \omega_T^\ell$.
Furthermore, we make the meaningful assumption that no patch $\omega$ coincides with the entire domain $\Omega$.

For its prototypical (global) basis functions $\{\varphi_j 
\with j =1,\dots,J \}$ associated to the element~$T$, the SLOD uses the following ansatz 
\begin{equation*}
	\varphi_j \coloneqq \mathcal A^{-1}g_j
\end{equation*}
with source terms $g_j\in \Pp(\omega)$ to be determined subsequently. We obtain a localized approximation $\psi_j
\in H^1_0(\omega)$ 
of the basis function $\varphi_j$ by computing its Galerkin projection onto the local subspace $H^1_0(\omega)$, i.e., $\psi_j
\in H^1_0(\omega)$ satisfies
\begin{equation}
	\label{e:patchproblem}
	a_\omega(\psi_j,v) = \tspf{g_j}{v}{\omega},\qquad\text{for all }v\in H^1_0(\omega).
\end{equation}

For the choice of $g_j$, we recall some notation and results on traces of $H^1(\omega)$-functions (see, e.g., \cite{LiM72a} for details). Denoting $U \coloneqq \V\lvert_\omega \subset H^1(\omega)$, we introduce the trace operator on $\omega$ restricted to $U$ as 
\begin{equation*}
	\operatorname{tr}\colon U \rightarrow X \coloneqq \mathrm{range}\,\operatorname{tr}\subset H^{1/2}(\partial  \omega).
\end{equation*}
An example of a continuous right-inverse of $\operatorname{tr}$ is the $\mathcal A$-harmonic extension, henceforth denoted by~$\operatorname{tr}^{-1}$. Given $w \in X$, it satisfies $\operatorname{tr} \operatorname{tr}^{-1}w= w$ and
\begin{equation}
	\label{eq:harmext}
	a_\omega(\operatorname{tr}^{-1}w,v)= 0,\qquad\text{for all }v \in H^1_0( \omega).
\end{equation}
Using definitions \cref{e:patchproblem,eq:harmext}, and that $(v - \operatorname{tr}^{-1}\operatorname{tr} v) \in H^1_0(\omega)$, it holds
\begin{equation*}
	a(\psi_j ,v) = a_\omega(\psi_j ,v) = a_\omega(\psi_j ,v-\operatorname{tr}^{-1} \operatorname{tr} v) = (g_j,v- \operatorname{tr}^{-1} \operatorname{tr} v)_{L^2(\omega)}.
\end{equation*}
This result yields, together with the definition of $\varphi_j$ and the  local support of $g_j$, the following key observation 
\begin{equation}
	\label{e:nd1}
	a(\varphi_j-\psi_j ,v) = (g_j,v)_{L^2(\omega)}-a_\omega(\psi_j ,v)= (g_j, \operatorname{tr}^{-1} \operatorname{tr}\, v)_{L^2(\omega)},\qquad v \in H^1_0(\Omega).
\end{equation}
Hence, we can rephrase the smallness of the localization error as the (almost) $L^2(\omega)$-orthogonality of $g_j$ to the space
\begin{equation}\label{e:Y}
	Y\coloneqq \operatorname{tr}^{-1}X\subset U
\end{equation}
of $\mathcal A$-harmonic functions on $\omega$ (which satisfy the homogeneous Dirichlet boundary condition on $\partial \Omega\cap\partial\omega$).
Since the restricted $L^2$-projection $\PiH\lvert_Y$ has a finite rank of dimension less or equal to $K \coloneqq J\cdot  \#(\T_{H}\cap \omega)$, there exists a singular value decomposition (SVD) such that
\begin{equation}
	\label{eq:op}
	\PiH\lvert_Y\,v = \sum_{k = 1}^{K} \sigma_k (v,w_k)_{H^1(\omega)}\,q_k
\end{equation}
where $\sigma_1\geq \dots \geq \sigma_{K}\geq 0$ denote the singular values, $\{q_1,\ldots,q_K\}$ the $L^2(\omega)$-orthonormal left singular vectors, and $\{w_1,\ldots,w_K\}$ the $H^1(\omega)$-orthonormal right singular vectors.

We choose the source terms $g_j$ as the left singular vectors corresponding to the $J$ smallest singular values, i.e.,
\begin{equation}
\label{eq:choiceg}
	g_j \coloneqq q_{K-J+j},\qquad j =1,\dots,J.
\end{equation}
This yields
\begin{equation*}
	\sup_{v\in Y\colon  \|v\|_{H^1(\omega)} = 1} {\tspf{g_j}{v}{\omega}} \leq \sigma_{K-J+1},\qquad j =1,\dots,J.
\end{equation*}
which follows directly from the properties of the SVD. We define the quantity $\sigma$ measuring the (quasi-)orthogonality between the $g_j$ and $Y$ as
\begin{equation}
	\sigma(H,\ell) \coloneqq \max_{T\in\TH}\sigma_{T}(H,\ell),\qquad
	\sigma_{T}(H,\ell)\coloneqq\sigma_{K-J+1}. \label{eq:sigma}
\end{equation}

The quantity $\sigma$ is crucial for the error analysis of the SLOD as it determines the localization error. Note that the dependence of $\sigma$ on the (fixed) polynomial degree $p$ is not made explicit in \cref{eq:sigma}. The following remark deals with the decay of $\sigma$ with respect to the oversampling parameter $\ell$ and, for the sake of curiosity, also the polynomial degree $p$.

\begin{remark}[Decay of $\sigma$]
	\label{r:decaysigma}
	In \cite{HaPe21b}, it has been numerically observed and conjectured that, for $p = 0$, the quantity $\sigma$ decays super-exponentially as $\ell$ is increased, cf.~\cref{fig:sv}~(left). A similar decay in $\ell$ can also be observed for $p>0$. In accordance with \cite{HaPe21b,Freese-Hauck-Peterseim}, we state the following conjecture: there exists $C_\sigma>0$ depending algebraically on $H, \ell$ and $C>0$ independent of $H,\ell$ such that 
	\begin{equation}
	\label{eq:conjdec}
	\sigma(H,\ell) \leq C_\sigma(H,\ell) \exp\big(-C \ell^{\frac{d-1}{d}}\big).
	\end{equation}

	Conversely, for fixed $\ell$, a rapid decay of $\sigma$ in $p$ can be observed, cf.~\cref{fig:sv}~(right). The low level of magnitude of the last singular values may suggest that the respective source terms correspond to fully local basis functions. However, due to the low levels or singular values, this is difficult to verify numerically. 
	\begin{figure}[h]
		\input{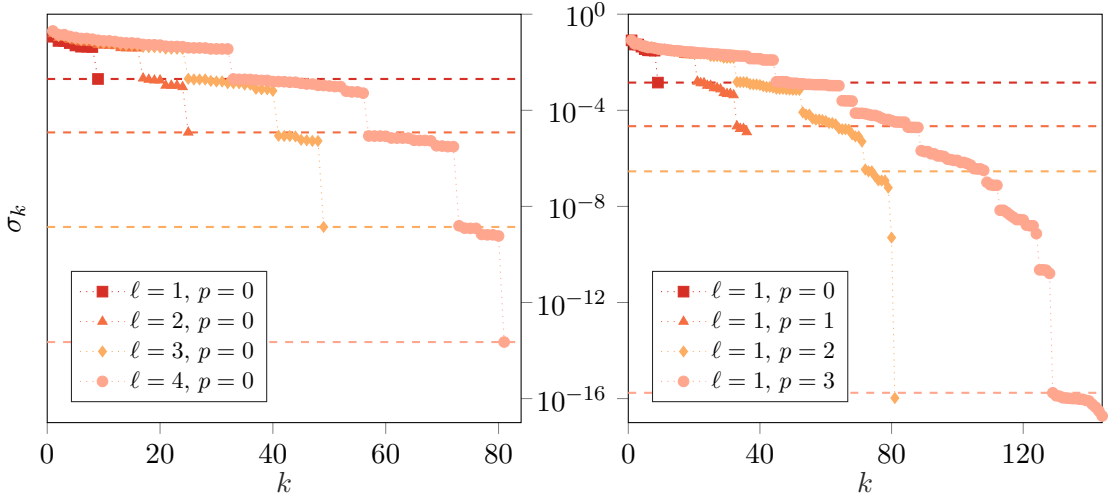}
		\caption{\small Singular values $\sigma_k$ of operator $\PiH\lvert_Y$ defined in \cref{eq:op} for an interior patch, for different pairs of $\ell,$ $p$. The singular values $\sigma_{K-J+1}$ relevant for~\cref{eq:sigma} are marked by dashed horizontal lines.}
		\label{fig:sv}
	\end{figure}
\end{remark}

Choosing for each patch the basis functions corresponding to the $J$ source terms~\cref{eq:choiceg}, we obtain the SLOD ansatz space
\begin{equation}
	\label{eq:Vslod}
	\Vslod \coloneqq \linh\{\pl\with \elem \in \TH, j = 1,\dots,J\}.
\end{equation}
The Galerkin SLOD solution $\vslod \in \Vslod$ then satisfies
\begin{equation} \label{eq:SLOD}
	a(\vslod,w) = \tspf{f}{w}{\Omega},\qquad \text{for all }w \in \Vslod.
\end{equation}

Note that a reasonable SLOD approximation requires a stable choice of the basis functions in~\cref{eq:Vslod}. However, for large oversampling parameters and patches intersecting the boundary, the choice \cref{eq:choiceg} may be insufficient, and a special treatment is required. In \cite[App.~B]{HaPe21b}, such an algorithm is proposed, curing possible stability and uniqueness issues in practice. Since we still cannot guarantee stability in an a priori manner, we assume that the source terms corresponding to the basis function of~\cref{eq:Vslod} form a Riesz basis of $\P(\TH)$, which is formulated in the following.

\begin{assumption}[Riesz stability]
	\label{a:riesz}
	The set
	\begin{equation*}
		\{\g\with \elem \in \TH, j = 1,\dots,J\}
	\end{equation*}
	is a Riesz basis of $\Pp(\TH)$, i.e., there is $C_{\mathrm{r}}(H,\ell)>0$ depending polynomially on $H$ and $\ell$ such that, for all $(c_{T,j})_{T \in \TH,j=1,\dots,J}$,
	\begin{equation}
		\label{e:Riesz}
		C^{-1}_{\mathrm{r}}(H,\ell)\sum_{\substack{\elem \in \mathcal T_H\\j=1,\dots,J}} c_{T,j}^2  \leq \bigg\lVert \sum_{\substack{\elem \in \mathcal T_H\\j=1,\dots,J}} c_{T,j} \g\biggl\rVert_{L^2(\Omega)}^2 \;\leq C_{\mathrm{r}}(H,\ell) \sum_{\substack{\elem \in \mathcal T_H\\j=1,\dots,J}} c_{T,j}^2.
	\end{equation}
\end{assumption}

\subsection{A posteriori error bound using SLOD}

We provide an a posteriori error analysis of the proposed method based on SLOD techniques. Conceptually, it is similar to the one for the SLOD, cf.~\cite[Thm.~6.1]{HaPe21b}, but it additionally includes the local optimal approximation error $\dn $ defined in \cref{e:kolmogorov}.

\begin{theorem}[A posteriori error bound]
	\label{t:aposteriori}
	Let \cref{a:riesz} be satisfied and let~$u$ and~$\vpum$ denote the solutions to \cref{eq:sol,eq:PUMSLOD}, respectively. Then, there exists a constant $\Ceo>0$ independent of $H,\ell,$ and $n$, such that, for any $f \in H^k(\TH)$, $k \in \mathbb N_0$,
	\begin{align*}
		\vnormof{u-\vpum}\leq \Ceo\Big( H^{s+1}|f|_{H^s(\TH)} + \ell^{d+1} C_\mathrm{r}^{1/2}(H,\ell) \big(\sigma(H,\ell) + H\dn(H,\ell)\big)\tnormf{f}{\Omega}\Big)
	\end{align*}
	with $s \coloneqq  \min\{k,p+1\}$ and the notation $H^0(\TH) \coloneqq  L^2(\Omega)$ and $\abs{\cdot}_{H^0(\TH)}\coloneqq  \|\cdot\|_{L^2(\Omega)}$. 
\end{theorem}
\begin{proof}
	The application of Céa's Lemma yields, for arbitrary $v \in \Vpum$,
	\begin{equation}
		\label{e:cea}
		\vnormof{u-\vpum} \leq \vnormof{u-v}.
	\end{equation}
	
	Let the solution to the (higher-order) collocation variant of the SLOD, cf.~\cite[Rem.~5.1]{HaPe21b}, with oversampling parameter $\ellp \coloneqq \lfloor\ell/2\rfloor$ be denoted by $\vslodp$. Given a source term $f \in L^2(\Omega)$, its solution is obtained by the following linear combination of localized basis functions $\plp$:
	\begin{equation*}
		\vslodp = \sum_{\substack{\elem \in \mathcal T_H\\j=1,\dots,J}} c_{T,j} \plp
	\end{equation*}
	with coefficients $c_{T,j}$ that are uniquely defined by
	\begin{equation}
		\sum_{\substack{\elem \in \mathcal T_H\\j=1,\dots,J}}c_{T,j} \gp = \Pi_H f. \label{eq:Pif}
	\end{equation}
	Adding and subtracting~$\vslodp$ in \cref{e:cea} and employing the triangle inequality yields
	\begin{equation}
		\label{e:trianglei}
		\vnormof{u-\vpum}  \leq \vnormf{u-\vslodp}{\Omega} + \vnormof{\vslodp-v}.
	\end{equation}
	
	The first term is the error of the collocation variant of the SLOD. It can be bounded using a higher-order version of \cite[Thm.~6.1]{HaPe21b} stating the existence of a constant $\Cslod>0$ independent of $H$ and $\ellp$ such that
	\begin{equation*}
		\vnormf{u-\vslodp}{\Omega}\leq \Cslod\big( H^{s+1}|f|_{H^s(\TH)}  +  C_{\mathrm{r}}^{1/2}(H,\ellp){\ellp}^{d/2}\sigma(H,\ellp)\|f\|_{L^2(\Omega)}\big).
	\end{equation*}
	
	For the second term in \cref{e:trianglei}, we choose $v \in \Vpum$ as sum of functions $v_z^n \in \Vpumz$ to be specified later, i.e., 
	\begin{equation*}
	\label{eq:v}
	v = \sum_{z \in \NH} \Lambda_z v_z^n.
	\end{equation*}
	Using the partition of unity property of the hat-functions $\sum_{z \in \NH}\Lambda_z \equiv 1$, we obtain
	\begin{align*}
		\vnormof{\vslodp -v}^2 = \Big\|\sum_{z \in \mathcal N_H} \Lambda_z(\vslodp - v_z^n)\Big\|_{a,\Omega}^2 \leq \Col \sum_{z \in \mathcal N_H}\vnormf{\Lambda_z(\vslodp - v_z^n)}{\omega_z}^2
	\end{align*}
	with $\Col$ defined in \cref{eq:overlap} denoting the maximal number of overlapping $\omega_z$. For any $z \in \NH$, we can locally on $\omega_z$ replace $\vslodp$ by $\vslodpz$ defined by
	\begin{equation*}
	\label{eq:slodsol4exch}
	\vslodpz \coloneqq \sum_{\substack{\elem \subset \omega_z^{\ellp}\\j=1,\dots,J}}c_{T,j}\plp \in H^1_0(\omega_z^\ell)
	\end{equation*}
	with the basis functions $\plp \coloneqq \mathcal A^{-1}_{\patchTm} \gp$. As approximation to $\vslodpz$, we use
	\begin{align*}
	v_z  \coloneqq \sum_{\substack{\elem \subset \omega_z^{\ellp}\\j=1,\dots,J}} c_{T,j} \plph \in \Vpumzi
	\end{align*}
	with $\plph \coloneqq \mathcal A^{-1}_{\omega_z^\ell} \gp \in \Vpumzi$ being approximations of $\plp$. We choose $v_z^n \in \Vpumz$ as the not necessarily unique element minimizing $\vnormf{\Lambda_z v_z -\Lambda_z v_z^n}{\omega_z^\ell}$.	Abbreviating 
	\begin{equation*}
	e_z^1 \coloneqq \vslodpz - v_z,\qquad e_z^2\coloneqq v_z - v_z^n,
	\end{equation*}
	and performing the above-mentioned local replacement, we obtain
	\begin{align}
		\label{e:esttriy}
		\vnormf{\Lambda_z(\vslodp - v_z^n)}{\omega_z}^2  
		 \leq 2\big( \vnormf{\Lambda_z e_z^1}{\omega_z}^2+\vnormf{\Lambda_z e_z^2}{\omega_z}^2\big),
	\end{align}
	where we add and subtract $v_z$ and employ the triangle inequality. 
	Using the product rule and the bound \cref{eq:boundhatfun}, we get for the first term
	\begin{align}
		\label{e:prodrule}
		\begin{split}
			& \vnormf{\Lambda_z e_z^1}{\omega_z}^2
			\leq   2\beta\big(\Cla^2 H^{-2} \tnormf{e_z^1}{\omega_z}^2 + \tnormf{\nabla e_z^1}{\omega_z}^2\big).
		\end{split}
	\end{align}
	Noting that by $e_z^1 \in H^1_0(\omega_z^\ell)$, we obtain for the first term in \cref{e:prodrule}
	\begin{align}
		\label{eq:l2toH1}
		\tnormf{e_z^1}{\omega_z}^2 \leq \tnormf{e_z^1}{\omega_z^{\ell}}^2 \leq \Cp^2 \pi^{-2}\ell^2 H^2 \tnormf{\nabla e_z^1}{\omega_z^\ell}^2
	\end{align}
	using Friedrichs' inequality on $\omega_z^\ell$ with $\operatorname{diam}(\omega_z^\ell) \leq \Cp \ell H$,  $\Cp>0$.
	For the second term in \cref{e:prodrule}, we infer the trivial estimate
	\begin{align*}
		\tnormf{\nabla e_z^1}{\omega_z}^2 \leq \tnormf{\nabla e_z^1}{\omega_z^\ell}^2,
	\end{align*}
	which implies that, in order to bound \cref{e:prodrule}, it suffices to estimate $\tnormf{\nabla e_z^1}{\omega_z^\ell}$. Using the continuity estimate
	\begin{equation*}
		\|\operatorname{tr}_{\patchTm}^{-1}\operatorname{tr}_{\patchTm} v\|_{H^1(\omega_T^m)} \leq C_\mathrm{tr} \|v\|_{H^1(\omega_T^m)},\qquad v \in \V\lvert_{\omega_T^m}
	\end{equation*}
	with a constant $C_\mathrm{tr}>0$ independent of $H$, $\ell$ from the proof of \cite[Thm.~6.1]{HaPe21b}, one can show that
	\begin{equation}
		\tspf{\gp}{\operatorname{tr}_{\patchTm}^{-1}\operatorname{tr}_{\patchTm} e_z^1}{\patchTm} \leq C_\mathrm{tr} \sigma(H,\ellp) \|e_z^1\|_{H^1(\patchTm)}.
		\label{eq:quasiorthogonality}
	\end{equation}
	By \cref{e:nd1,eq:quasiorthogonality},  as well as the discrete Cauchy--Schwarz inequality, the finite overlap of the patches $\patchTm$ and Friedrichs' inequality, we get 
	\begin{align*}
		\alpha \tnormf{\nabla e_z^1}{\omega_z^\ell}^2 &\leq \sum_{\substack{\elem \subset \omega_z^{\ellp}\\j=1,\dots,J}} c_{T,j} a(\plp-\plph,e_z^1) = -\sum_{\substack{\elem \subset \omega_z^{\ellp}\\j=1,\dots,J}} c_{T,j} \tspf{\gp}{\operatorname{tr}_{\patchTm}^{-1}\operatorname{tr}_{\patchTm} e_z^1}{\patchTm}\\
		&\leq C_\mathrm{tr}\sigma(H,\ellp) \sum_{\substack{\elem \subset \omega_z^{\ellp}\\j=1,\dots,J}} \abs{c_{T,j}} \|e_z^1\|_{H^1(\patchTm)}\\
		&\leq C C_\mathrm{tr}  (\Cp \ell H \pi^{-1} + 1)\sigma(H,\ellp) {\ellp}^{d/2}J^{1/2}\sqrt{\sum_{\substack{\elem \subset \omega_z^{\ellp}\\j=1,\dots,J}} c_{T,j}^2}\; \tnormf{\nabla e_z^1}{\omega_z^\ell}
	\end{align*}
	with constant $C>0$ reflecting the overlap of the patches~$\omega_T^\ellp$. Using \cref{eq:l2toH1}, this yields an estimate for \cref{e:prodrule} and consequently bounds the first term in \cref{e:esttriy}.
	
	For the second expression in \cref{e:esttriy}, using the definition of the Kolmogorov $n$-width~\cref{Kolmogorovnw}, Friedrichs' inequality on the patch $\omega_z^\ell$, the discrete Cauchy--Schwarz inequality, and that the $\gp$ are $L^2$-normalized, we obtain
	\begin{align*}
		\vnormf{\Lambda_ze_z^2}{\omega_z} &\leq \dn(H,\ell)\vnormf{v_z}{\omega_z^\ell} = \dn(H,\ell)\bigg\|\mathcal A_{\omega_z^\ell}^{-1}\sum_{\substack{\elem \subset \omega_z^{\ellp}\\j=1,\dots,J}} c_{T,j} \gp\bigg\|_{a,\omega_z^\ell}\\
		&\leq  \Cp\dn(H,\ell)\alpha^{-1/2} \pi^{-1} \ell H \bigg\|\sum_{\substack{\elem \subset \omega_z^{\ellp}\\j=1,\dots,J}}c_{T,j} \gp\bigg\|_{L^2(\omega_z^\ell)}\\
		&\leq C \Cp\dn(H,\ell)\alpha^{-1/2} \pi^{-1} \ell H  \ellp^{d/2}J^{1/2} \sqrt{\sum_{\substack{\elem \subset \omega_z^{\ellp}\\j=1,\dots,J}}c_{T,j}^2}
	\end{align*}
	with constant $C>0$ appearing in the bound $C^2 \ellp^dJ$ of the number of terms in the above sum. Consequently, we conclude the estimate for the second term in~\cref{e:esttriy}.
	
	The assertion can be finalized by combining all estimates utilizing
	\begin{align*}
		\sum_{z \in \mathcal N_H} \sum_{\substack{\elem \subset \omega_z^{\ellp}\\j=1,\dots,J}}c_{T,j}^2  \leq C m^d \sum_{\substack{\elem \in \TH \\j=1,\dots,J}} c_{T,j}^2 \leq C C_\mathrm{r}(H,\ell) m^d \tnorm{f}{\Omega}^2
	\end{align*}
with constant $C>0$ reflecting the overlap of the patches $\omega_z^m$. Here, we used \cref{a:riesz,eq:Pif,e:stabL2}. Finally, for the sake of readability, we substitute $\ellp = \lfloor \tfrac{\ell}{2}\rfloor$ by $\ell$, which may introduce additional constants that change the decay rate of $\sigma$ by some constant factor.
\end{proof}

\subsection{Local approximation error bound using SLOD} 
The error estimate from \cref{t:aposteriori} incorporates the spectral approximation error $\dn$ defined in \cref{e:kolmogorov}. Subsequently, we derive a bound for $\dn$ utilizing the SLOD for constructing bases of the local spaces $\Vpumzi$ defined in \cref{eq:V_H_ell_z}. For this purpose, we fix a node $z \in \NH$ and treat $\omega_z^ \ell$ as the whole domain. For the oversampling parameters $m = 1,\dots,\ell-1$, we denote the basis functions with a tilde to emphasize that, in general, they do not coincide with their global counterparts~$\plp$,~i.e., 
\begin{equation}
\label{eq:philoc}
	\{\ppt \coloneqq \mathcal A_{\omega_z^\ell}^{-1} \gpt\with T \subset \omega_z^\ell,\,j = 1,\dots,J\}
\end{equation}
with source terms $\gpt \in \P(\tilde \omega_T^m)$, where $\tilde \omega_T^m \coloneqq \omega_T^m \cap \omega_z^\ell$. We denote the corresponding localized basis by
\begin{equation}
\label{eq:psiloc}
	\{\plpt \coloneqq \mathcal A_{\patchTmp}^{-1} \gpt\with T \subset \omega_z^\ell,\,j = 1,\dots,J\}.
\end{equation}

Similarly, as in the full domain setting, we need to measure, for all $T \subset \omega_z^\ell$, the \mbox{(quasi-)} orthogonality of the source terms $\{\gpt\with j = 1,\dots,J\}$ on the corresponding space of $\mathcal A$-harmonic functions
\begin{equation}
	 \operatorname{tr}^{-1}_{\patchTmp}\operatorname{tr}_{\patchTmp}\big(H^1_0(\omega_z^\ell)\lvert_{\patchTmp}\big).\label{e:Ypatchlocal}
\end{equation}
Similarly to \cref{eq:op}, for a node $z \in \NH$ and element $T \subset \omega_z^\ell$, we denote the singular values of $\PiH$ restricted to~\cref{e:Ypatchlocal} by $\tilde \sigma_1\geq \tilde \sigma_2\geq  \dots\geq \tilde \sigma_K\geq 0$. Analogously to \cref{eq:sigma}, we define 
\begin{equation*}
	\tilde \sigma(H,\ell,m) \coloneqq \max_{z\in\NH}\max_{T\in\TH} \tilde \sigma_{z,T},\qquad\tilde \sigma_{z,T} \coloneqq \tilde \sigma_{K-J+1}.
\end{equation*}

The quantity $\tilde \sigma$ is strongly related to its counterpart $\sigma$ from \cref{eq:sigma} and, in numerical experiments, exhibits the same qualitative behavior as described in \cref{r:decaysigma}. A local variant of \cref{a:riesz} is required to ensure the stability of the local basis. 

\begin{assumption}[Local Riesz stability]
	\label{a:rieszloc}
	For all patches $\omega_z^\ell$, the set
	\begin{equation*}
		\{\gpt \with \elem \subset \omega_z^\ell, j = 1,\dots,J\}
	\end{equation*}
	is a Riesz basis of $\Pp(\omega_z^\ell)$, i.e., there is $\tilde C_{\mathrm{r}}(H,\ell,\ellp)>0$ depending polynomially on $H,\ell,$ and $\ellp$ such that, for all $z \in \NH$ and all $(c_{T,j})_{\elem\subset \omega_z^\ell,j=1,\dots,J}$,
	\begin{equation}
		\label{e:Rieszloc}
		\tilde C^{-1}_{\mathrm{r}}(H,\ell,\ellp)\sum_{\substack{\elem \subset \omega_z^\ell\\j=1,\dots,J}} c_{T,j}^2  \leq \bigg\lVert \sum_{\substack{\elem \subset \omega_z^\ell\\j=1,\dots,J}} c_{T,j} \gpt\biggl\rVert_{L^2(\omega_z^\ell)}^2 \;\leq \tilde C_{\mathrm{r}}(H,\ell,\ellp) \sum_{\substack{\elem\subset \omega_z^\ell\\j=1,\dots,J}} c_{T,j}^2.
	\end{equation}
\end{assumption}

\begin{theorem}[Bound on $\dn$]\label{thm:bound_delta_SLOD}
	Let \cref{a:rieszloc} be satisfied. Then, there exists a constant $\Cdo>0$ independent of $H,\ell,$ and $m$ such that, for $\ellp = 1,\dots,\ell-1$ 
	\begin{equation*}
		\dn(H,\ell)
		\leq  \Cdo\ell\, \ellp^{d/2}  H^{-1} \tilde C^{1/2}_\mathrm{r}(H,\ell,\ellp) \tilde\sigma(H,\ell,\ellp),
	\end{equation*}
	where $n \approx m^d$.
\end{theorem}

\begin{proof}
	Let us consider a fixed node $z \in \NH$ and oversampling parameter $\ell$.
	As approximation space $Q(n)$ of dimension $n \approx \ellp^d$, we choose
	\begin{equation*}
		Q(n) \coloneqq \linh\{\Lambda_z \plpt\with T \subset \omega_z^\ellp, j = 1,\dots,J\}
	\end{equation*} 
	with basis functions defined in \cref{eq:psiloc}. For the approximation of $v_z \in \Vpumzi$, we choose the element	$w_z \in Q(n)$ as 
	\begin{equation*}
		w_z = \Lambda_z \sum_{\substack{\elem \subset \omega_z^\ellp\\j=1,\dots,J}} c_{T,j}\plpt,
	\end{equation*}
	where the $c_{T,j}$ are the coefficients of the expansion of $v_z$ in terms of the basis functions~$\ppt$ defined in \cref{eq:philoc}. Note that, by \cref{a:rieszloc}, the coefficients $c_{T,j}$ are uniquely determined. Thus, we can estimate the spectral approximation error \cref{Kolmogorovnw} using $v_z$ and~$w_z$~as
	\begin{align*}
		\dnz(H,\ell) = \adjustlimits \inf_{Q(n) \subset H^1_0(\omega_z)} \sup_{v_z \, \in \, \Vpumzi} \inf_{w_z \in Q(n)} \frac{\vnormf{\Lambda_zv_z -w_z}{\omega_z}}{\vnormf{v_z}{\omega_z^\ell}} \leq \sup_{v_z \, \in \, \Vpumzi} \frac{\vnormf{\Lambda_zv_z -w_z}{\omega_z}}{\vnormf{v_z}{\omega_z^\ell}}.
	\end{align*}

	Denoting
	\begin{equation*}
		e_z \coloneqq \sum_{\substack{\elem \subset \omega_z^\ellp\\j=1,\dots,J}} c_{T,j}\big(\ppt-\plpt\big) \in H^1_0(\omega_z^\ell),
	\end{equation*}
	we obtain for the numerator, using the product rule, the triangle inequality, and \cref{eq:boundhatfun}:
	\begin{align*}
		\vnormf{\Lambda_zv_z -w_z}{\omega_z}\leq \beta^{1/2} \tnorm{\nabla (\Lambda_z e_z)}{\omega_z^\ell} \leq \beta^{1/2}\big(\Cla H^{-1} \tnorm{e_z}{\omega_z^\ell} + \tnorm{\nabla e_z}{\omega_z^\ell}\big).
	\end{align*}
	We apply Friedrichs' inequality on the patch $\omega_z^\ell$ using that $\operatorname{diam}(\omega_z^\ell)\leq \Cp\ell H$,  with a constant $\Cp>0$. Hence, we can bound the first term against the second term, i.e.,
	\begin{align*}
		\tnorm{e_z}{\omega_z^\ell} \leq \Cp H \ell \pi^{-1} \tnorm{\nabla e_z}{\omega_z^\ell}.
	\end{align*}
We adapt estimate \cref{eq:quasiorthogonality} to the local setting introducing a constant $\tilde C_\mathrm{tr}>0$. Using the finite overlap of the patches $\patchTmp$, the discrete Cauchy--Schwarz inequality and Friedrichs' inequality on $\omega_z^\ell$, we obtain
	\begin{align*}
		\alpha\tnormf{\nabla e_z }{\omega_z^\ell}^2 &\leq \sum_{\substack{\elem \subset \omega_z^\ellp\\j=1,\dots,J}} c_{T,j} a(\ppt-\plpt,e_z) \leq  \tilde C_\mathrm{tr} \tilde\sigma(H,\ell,\ellp) \sum_{\substack{\elem \subset \omega_z^{\ellp}\\j=1,\dots,J}}c_{T,j} \|e_z\|_{H^1(\patchTmp)}\\
		&\leq C \tilde C_\mathrm{tr} (\Cp \pi^{-1}H\ell + 1) \ellp^{d/2} \tilde\sigma(H,\ell,\ellp) \sqrt{\sum_{\substack{\elem \subset \omega_z^\ellp\\j=1,\dots,J}} c_{T,j}^2}\|\nabla e_z\|_{L^2(\omega_z^\ell)},
	\end{align*}
	where $C>0$ reflects the overlap of the patches $\patchTmp$. 
	
	Adding the remaining coefficients $c_{T,j}$ from the expansion of $v_z$ in terms of the $\ppt$ and using \cref{a:rieszloc}, we get
	\begin{align*}
		\tilde C_\mathrm{r}^{-1}(H,\ell,\ellp)\sum_{\substack{\elem \subset \omega_z^\ell\\j=1,\dots,J}}c_{T,j}^2 &\leq  \bigg\lVert \sum_{\substack{\elem \subset \omega_z^\ell\\j=1,\dots,J}} c_{T,j} \gpt \biggl\rVert_{L^2(\omega_z^\ell)}^2
		\leq \Cbo^2  H^{-2} \bigg\lVert \sum_{\substack{\elem \subset \omega_z^\ell\\j=1,\dots,J}} c_{T,j} \gpt \biggl\rVert_{H^{-1}(\omega_z^\ell)}^2\\&
		\leq \Cbo^2  \beta H^{-2} \vnorm{v_z}{\omega_z^\ell}^2.
	\end{align*}
	Here, we also employed that, by \cref{eq:BHL2proj,eq:bubbleopest}, we have, for all $q \in \P(\omega_z^\ell)$,
	\begin{align*}
		\tnormf{q}{\omega_z^\ell} &= \sup_{v \, \in \, H^1_0(\omega_z^\ell)} \frac{\tspf{q}{v}{\omega_z^\ell}}{\tnormf{v}{\omega_z^\ell}} \leq \Cbo H^{-1}\sup_{v \, \in \, H^1_0(\omega_z^\ell)} \frac{\tspf{q}{\BH v}{\omega_z^\ell}}{\tnormf{\nabla \BH v}{\omega_z^\ell}}\\
		&\leq \Cbo H^{-1} \|q\|_{H^{-1}(\omega_z^\ell)}.
	\end{align*}
	Combining the estimate yields the assertion.
\end{proof}

\begin{remark}[Choice of parameters]
	\label{rem:bound_delta_SLOD}
	This remark specifies how to choose the oversampling parameter~$\ell$ and the number of local functions~$n$ in order to preserve the optimal order of convergence of $s+1$ in \cref{t:aposteriori}.
	For $\ell$, the super-exponential decay~\cref{eq:conjdec} implies that it needs to be chosen of order $\abs{\log H}^{(d-1)/d}$. Using \cref{thm:bound_delta_SLOD} and that $\tilde \sigma$ has similar decay properties as $\sigma$, we obtain that $n$ needs to be chosen of order $\abs{\log H}^{d-1}$.
	Note that these choices require the validity of~\cref{eq:conjdec} and \cref{a:riesz,a:rieszloc}. For a (pessimistic) estimate which is valid without additional assumptions, we refer to \cref{rem:bound_delta_LOD} below.
\end{remark}

\section{(Pessimistic) a priori error analysis} \label{sec:a priori error analysis}

This section presents an a priori error analysis of the proposed method, which is based on the LOD framework, cf.~\cite{MaP14,HeP13,Peterseim2021}. Note that the exponential localization properties of the LOD cannot match the practically observed super-exponential localization properties of the SLOD, cf.~\cref{r:decaysigma}. Nevertheless, the LOD construction has the decisive advantage that the basis stability is guaranteed by construction and that the exponential localization can be rigorously proved. This enables an a priori analysis without assumptions on the stability of the SLOD basis and without conjectures on the decay of singular values,  cf.~\cref{a:riesz,a:rieszloc,r:decaysigma}.

\subsection{Higher-order LOD}
We briefly introduce a higher-order version of the LOD similar to the constructions in \cite{Mai20ppt,DHM22}. 
The LOD constructs its problem-adapted basis functions by adding fine-scale information to coarse-scale finite element functions. We define the space of fine-scale functions~as 
\begin{equation}
	\W \coloneqq \mathrm{kernel}\,\Pi_H.
	\label{eq:W}
\end{equation}

The step of adding fine-scale information is called correction and utilizes the so-called correction operator $\C\colon \V \to \W$ defined as the $a$-orthogonal projection onto~$\W$, i.e.,
\begin{equation*}
	a(\C v,w) = a(v,w),\qquad \text{for all }w\in \W.
\end{equation*}
We split up the correction operator into a sum of element correction operators, i.e., $\C = \sum_{T \in \TH} \Ct$ with element correction operators $\Ct\colon \V \to \W$ defined by
\begin{equation*}
	a(\Ct v,w) = a_T(v,w),\qquad \text{for all }w\in \W.
\end{equation*}
For any $v \in \V$, the correction $\Ct v$ decays exponentially fast away from its associated element $T \in \TH$, cf.~\cite[Lemma~5.1]{DHM22}, which motivates a localization. For this purpose, we substitute the global space $\W$ by localized counterparts $\W_T^\ell \coloneqq \W(\omega_T^\ell)$, where, for a subdomain $\omega \subset \Omega$, we use the definition
\begin{equation}
\label{eq:defrestW}
	\W(\omega) \coloneqq \{w \in \W\with \supp(w)\subset \omega\}.
\end{equation}
The localized element correction operator $\Clt\colon \V \to \W_T^{\ell}$ is then defined by
\begin{equation*}
	a(\Clt v,w) = a_T(v,w),\qquad \text{for all }w\in \W_T^\ell.
\end{equation*}
Similarly as for the correction operator $\C$ which can be decomposed into a sum of element correction operators, we define the localized correction operator by 
\begin{equation*}
	\Cl \coloneqq \sum_{T \in \TH} \Clt.
\end{equation*}

The ansatz space of the LOD is then obtained by adding (localized) corrections to the bubble functions $\{\bTj\with T\in \TH, j = 1,\dots,J\}$ defined in \cref{eq:bubble}, i.e., 
\begin{equation*}
	\Vlod \coloneqq \linh\{(1-\Cl) b_{T,j}\with T \in \TH, j=1,\dots,J\}.
\end{equation*} 
The Galerkin LOD approximation $\vlod \in \Vlod$ then satisfies
\begin{equation}
	\label{eq:LODsol}
	a(\vlod,w) = \tspf{f}{w}{\Omega},\qquad \text{for all }w \in \Vlod.
\end{equation}

\subsection{A priori error bound using LOD}

Using LOD techniques, we can prove the following a priori error estimate for the proposed multi-scale method. 
Numerical experiments show that this estimate is tentatively pessimistic. Nevertheless, compared to \cref{t:aposteriori}, it has the crucial advantage that it does not rely on additional assumptions or conjectures.

\begin{theorem}[A priori error bound]
	\label{t:apriori}
	Let $u$ and $\vpum$ denote the solutions of \cref{eq:sol,eq:PUMSLOD}, respectively. There exist constants $\Cet,\Cd>0$ independent of $H,\ell,$ and $n$ such that, for any $f \in H^k(\TH)$, $k \in \mathbb N_0$
	\begin{align*}
		\vnormof{u-\vpum} &\leq \Cet\Big( H^{s+1}|f|_{H^s(\TH)} + \ell^{d/2}H^{-1} \big(\ell^{d/2}\exp(-\Cd \ell) + \dn(H,\ell)\big)\tnormf{f}{\Omega}\Big).
	\end{align*}
	with $s \coloneqq  \min\{k,p+1\}$.
\end{theorem}

\begin{proof}
	We apply Céa's Lemma, which yields, for arbitrary $v \in \Vpum$,
	\begin{equation}
	\label{eq:cealod}
		\vnormof{u-\vpum} \leq \vnormof{u-v}.
	\end{equation}
	For the oversampling parameter $\ellp \coloneqq \lfloor\ell/2\rfloor$, we define the approximation
	\begin{equation*}
		\vlodp \coloneqq (1-\C^\ellp)\BH u
	\end{equation*}
	which is not the Galerkin LOD solution \cref{eq:LODsol} but has the same approximation properties, cf.~proof of \cite[Thm.~4.4]{Mai20ppt}. Adding and subtracting $\vlodp$ in~\eqref{eq:cealod} and using the triangle inequality yields
	\begin{equation}
		\label{e:trianglei2}
		\vnormof{u-\vpum}  \leq \vnormf{u-\vlodp}{\Omega} + \vnormof{\vlodp-v}.
	\end{equation}
	
	Using the above-mentioned approximation properties of $\vlodp$, we obtain the following estimate for the first term in \cref{e:trianglei2}
	\begin{equation*}
		\vnormf{u-\vlodp}{\Omega}\leq \Clod\big( H^{s+1}|f|_{H^s(\TH)}  +  \ellp^{d/2} H^{-1}\exp(-\Cd\ellp)\|f\|_{L^2(\Omega)}\big).
	\end{equation*}
	with constants $\Clod,\Cd>0$ independent of $H,\ell,$ and $m$.
	
	For the second term in \cref{e:trianglei2}, we choose $v \in \Vpum$ as sum of functions $v_z^n \in \Vpumz$ to be specified later, i.e., 
	\begin{equation*}
	\label{eq:vlod}
	v = \sum_{z \in \NH} \Lambda_z v_z^n.
	\end{equation*}
	Using the partition of unity property of the hat-functions $\sum_{z \in \NH}\Lambda_z \equiv 1$, we obtain for the second term in~\cref{e:trianglei2},
	\begin{align*}
	\vnormof{\vlodp -v}^2 = \Big\|\sum_{z \in \mathcal N_H} \Lambda_z(\vlodp - v_z^n)\Big\|_{a,\Omega}^2 \leq \Col \sum_{z \in \mathcal N_H}\vnormf{\Lambda_z(\vlodp - v_z^n)}{\omega_z}^2
	\end{align*}
	with $\Col$ defined in \cref{eq:overlap}. For any $z \in \NH$, we can, locally on $\omega_z$, substitute $\vlodp$ by
	\begin{equation*}
\vlodpz \coloneqq (1-\C^\ellp)(\BH u\lvert_{\omega_z^\ellp}).
\end{equation*}
Hence, we define $v_z^n\in  \Vpumz$ as the (not necessarily unique) elements minimizing the expression $\vnormf{\Lambda_z v_z -\Lambda_z v_z^n}{\omega_z^\ell}$, where
	\begin{align*}
		v_z \coloneqq (1-\Clzt)(\BH u\lvert_{\omega_z^\ellp})
	\end{align*}
	is an approximation to $\vlodpz$. We denote $\W_z^\ell\coloneqq \W(\omega_z^\ell)$ and define the above used correction operator $\Clzt\colon H^1_0(\Omega)\to \W_z^\ell$ by
	\begin{equation}
		a_{\omega_z^\ell}( \Clzt v,w) = a_{\omega_z^\ell}(v,w)\qquad\text{for all }w \in \W(\omega_z^\ell)
		\label{eq:Clz}.
	\end{equation}
	Note that it holds $v_z \in \Vpumzi$ which is a non-trivial observation, cf.~\cite[Rem~3.7~ii]{Peterseim2021}.
	
Abbreviating 
\begin{equation*}
e_z^1 \coloneqq \vlodpz - v_z,\qquad e_z^2\coloneqq v_z - v_z^n,
\end{equation*}
we obtain after performing the above-mentioned local substitution
	\begin{align}
\label{e:esttriylod}
\vnormf{\Lambda_z(\vlodp - v_z^n)}{\omega_z}^2  
\leq 2\big( \vnormf{\Lambda_z e_z^1}{\omega_z}^2+\vnormf{\Lambda_z e_z^2}{\omega_z}^2\big).
\end{align}
Following \eqref{e:prodrule}, in order to bound the first term in \cref{e:esttriylod}, it suffices to bound $\tnormf{e_z^1}{\omega_z}$ and~$\tnormf{\nabla e_z^1}{\omega_z}$. It holds
	\begin{equation*}
	e_z^1 = (\C^\ellp - \Clzt) (\BH u\lvert_{\omega_z^\ellp}) \in \W_z^\ell
	\end{equation*}
	which implies that $e_z^1$ has vanishing element averages. Thus, by Poincare's inequality
	\begin{align*}
		\tnormf{e_z^1}{\omega_z}^2 \leq 4\pi^{-2}H^2 \tnormf{\nabla e_z^1}{\omega_z}^2
	\end{align*}
	it is sufficient to estimate $\tnormf{\nabla e_z^2}{\omega_z}$ in order to obtain a bound for the first term in \eqref{e:esttriylod}. Given a function supported in ${\omega_z^\ellp}$ (e.g., $\BH u\lvert_{\omega_z^\ellp}$) the correction operator $\C^\ellp$ coincides with the localization of $\Clzt$ to $\ellp$-th order patches. Hence, we can apply the localization error estimate from the proof of \cite[Thm.~4.4]{Mai20ppt}, here in the oversampling parameter $\ellp$ and~\cref{eq:bubbleopest} to obtain
	\begin{align*}
		\tnormf{\nabla e_z^1}{\omega_z}^2 &\leq \tnormf{\nabla (\C^\ellp-\Clzt)(\BH u\lvert_{\omega_z^\ellp})}{\omega_z^\ell}^2\leq \Cloc^2 \ellp^{d}\exp(-\Cd\ellp)^2 \tnormf{\nabla (\BH u\lvert_{\omega_z^\ellp})}{\omega_z^\ell}^2\\
		&\leq \Cloc^2 \Cbo^2 H^{-2} \ellp^{d} \exp(-\Cd\ellp)^2 \tnormf{u}{\omega_z^\ellp}^2.
	\end{align*}
	with constants $\Cloc,\Cd>0$ independent of $H,\ell,$ and $m$. 
	
	For the second term in 	\cref{e:esttriylod}, we obtain by the definition of the Kolmogorov $n$-width, the stability of $(1-\Clzt)$ and \cref{eq:bubbleopest} that
	\begin{align*}
		\vnormf{\Lambda_z e_z^2}{\omega_z}^2 &\leq \dn^2(H,\ell,n)\vnormf{v_z}{\omega_z^\ell}^2
		\leq \beta \dn^2(H,\ell,n)\tnormf{\nabla (\BH u\lvert_{\omega_z^\ellp})}{\omega_z^\ell}^2\\*
		&\leq \Cbo^2 H^{-2}\beta^2 \dn^2(H,\ell,n) \tnormf{u}{\omega_z^\ellp}^2.
	\end{align*}
	Combining the previous estimates, using Friedrichs' inequality on $\Omega$ (recall that $\Omega$ is scaled to unit size), we get
	\begin{equation*}
		\tnormf{u}{\Omega} \leq \pi^{-1}\tnormf{\nabla u}{\Omega} \leq \pi^{-2}\alpha^{-1}\tnormf{f}{\Omega}.
	\end{equation*}
	Finally, we substitute $m = \lfloor \tfrac{\ell}{2}\rfloor$ by $\ell$, which introduces additional constants and changes the exponential decay rate $\Cd$ by a factor of two.
\end{proof}

\subsection{Local approximation error bound using LOD}

Subsequently, we derive an a priori bound for $\dn$, which is fully explicit in $H$ and $\ell$. This is the analog to \cref{thm:bound_delta_SLOD}, which does not rely on additional assumptions or conjectures. Numerical experiments show that this estimate is tentatively pessimistic. 

\begin{theorem}[Bound on $\dn$]
	\label{thm:bound_delta_LOD}
	There exists $\Cdt,\Cd>0$ independent of $H,\ell,$ and $m$ such that, for $\ellp = 1,\dots,\ell-1$ 
	\begin{equation*}
		\dn(H,\ell)
		\leq  \Cdt\ell\, \ellp^{d/2} \exp(-\Cd m),
	\end{equation*}
	where $n \approx m^d$.
\end{theorem}
\begin{proof}
	Let us consider a fixed node $z \in \NH$ and oversampling parameter $\ell$. By \cite[Rem.~3.7~ii]{Peterseim2021}, we can write any $v \in \Vpumzi$ as 
	\begin{align*}
		v = (1-\Clzt)\BH v
	\end{align*}
	with the correction operator $\Clzt$ defined in \cref{eq:Clz}. 
	We define the patch-local localized correction operator by
	$\tilde \C^\ellp \coloneqq \sum_{T \subset \omega_z^\ell} \tilde \C_T^\ellp$, where, denoting $\tilde \W_T^m \coloneqq \W( \omega_T^\ellp\cap \omega_z^\ell)$, the element correctors $\tilde \C_T^\ellp\colon H^1_0(\omega_z^\ell)\to \tilde \W_T^m$ are defined by
	\begin{equation*}
		a_{\omega_z^\ell}(\tilde \C_T^\ellp v,w) = a_T(v,w),\qquad\text{for all }w \in\tilde \W_T^m.
	\end{equation*}

	As approximation space $Q(n)$ of dimension $n \approx \ellp^d$, we choose
	\begin{equation*}
		Q(n) \coloneqq \linh\{\Lambda_z (1-\tilde \C^\ellp)b_{T,j}\with T \subset \omega_z^\ellp, j = 1,\dots,J\}
	\end{equation*} 
	and as approximation $w_z \in Q(n)$ of an element $v_z\in \Vpumzi$, we use
	\begin{equation*}
		w_z = \Lambda_z (1-\tilde \C^\ellp)(\BH v_z\lvert_{\omega_z^\ellp}).
	\end{equation*}
	
	Using the approximation space $Q(n)$ and the above defined choice of $w_z$, we can bound the Kolmogorov $n$-width as follows
	\begin{align*}
		\dnz(H,\ell) = \adjustlimits \inf_{Q(n) \subset H^1_0(\omega_z)} \sup_{v_z \, \in \, \Vpumzi} \inf_{w_z \in Q(n)} \frac{\vnormf{\Lambda_z v_z -w_z}{\omega_z}}{\vnormf{v_z}{\omega_z^\ell}} \leq \sup_{v_z \, \in \, \Vpumzi} \frac{\vnormf{\Lambda_z v_z -w_z}{\omega_z}}{\vnormf{v_z}{\omega_z^\ell}}.
	\end{align*}
	Abbreviating
	\begin{equation*}
	e_z \coloneqq (\tilde \C^\ellp-\Clzt) \BH v_z,
	\end{equation*}
	we can estimate the numerator using \cref{eq:boundhatfun} and 
	\begin{equation*}
		\Lambda_z(1-\tilde \C^\ellp)(\BH v_z\lvert_{\omega_z^\ellp}) = \Lambda_z(1-\tilde \C^\ellp)\BH v_z,
	\end{equation*}
	as 
	\begin{align*}
		\vnormf{\Lambda_zv_z -w_z}{\omega_z}& 
		\leq  \beta^{1/2}\tnormf{\nabla \Lambda_z e_z}{\omega_z}\leq \beta^{1/2}\big(\Cla H^{-1} \tnormf{e_z}{\omega_z} + \tnormf{\nabla e_z}{\omega_z}\big).
	\end{align*}
	It holds that $e_z \in \W_z^\ell$, which implies that $e_z^1$ has vanishing element averages. Thus, by Poincare's inequality
	\begin{equation*}
		\tnormf{e_z}{\omega_z} \leq 2\pi^{-1}H\tnormf{\nabla e_z}{\omega_z}.
	\end{equation*}
	Applying the localization error estimate from the proof of \cite[Tehorem 4.4]{Mai20ppt} to show that $\tilde \C^\ellp$ approximates $\Clzt$ exponentially and using \cref{eq:bubbleopest}, and Friedrichs' inequality on the patch $\omega_z^\ell$ with $\operatorname{diam}(\omega_z^\ell) \leq \Cp \ell H$, $\Cp >0$, we obtain 
	\begin{align*}
		\vnormf{\Lambda_zv_z -w_z}{\omega_z}&\leq\beta^{1/2} (2\Cla\pi^{-1}+1)\tnormf{\nabla e_z}{\omega_z^\ell}\\
		&\leq \Cloc m^{d/2} \beta^{1/2}(2\Cla \pi^{-1}+1)\exp(-\Cd m) \tnormf{\nabla\BH v}{\omega_z^\ell}\\
		&\leq \Cbo\Cloc\Cp \ell\, \ellp^{d/2} \beta^{1/2}(2\Cla \pi^{-1}+1)\alpha^{-1/2}\pi^{-1}\exp(-\Cd m) \vnormf{v_z}{\omega_z^\ell},
	\end{align*}
	with constants $\Cloc,\Cd>0$ independent of $H,\ell$, and $m$. The assertion follows immediately.
\end{proof}

\begin{remark}[Choice of parameters]
	\label{rem:bound_delta_LOD}
	This remark specifies how to choose the oversampling parameter~$\ell$ and the number of local functions~$n$ in order to guarantee the optimal order of convergence of $s+1$ in \cref{t:apriori}. By \cref{t:apriori}, $\ell$ needs to be chosen of order $\abs{\log H}$. Using \cref{thm:bound_delta_LOD}, we obtain that $n$ needs to be chosen of order $\abs{\log H}^{d}$.
	According to the experiments, these choices are pessimistic, cf. \cref{rem:bound_delta_SLOD}.
\end{remark}
\section{Implementation and numerical experiments}\label{sec:Implementation and numerical experiments}

In this section, we numerically investigate the proposed multi-scale method regarding the localization error, optimal convergence properties, high-contrast channeled coefficients, and higher-order polynomials using suitable benchmark problems. As a comparison, we use the SLOD from \cite{HaPe21b}, which we consider as state-of-the-art. We refer to \cite[Sec.~8]{HaPe21b} for a comparison of the SLOD to other multi-scale methods such as the LOD. For underlining the origin of the proposed method and its super-localization properties, cf.~\cref{t:aposteriori,thm:bound_delta_SLOD}, we subsequently refer to it as Super-Localized Generalized Finite Element Method~(SL-GFEM).

\subsection{Implementation}

For the practical implementation of the SL-GFEM, we need to perform a fine-scale discretization, i.e., we substitute the infinite-dimensional function space~$\V$ by the finite element space $\V_h \coloneqq \V \cap \P_1(\Th)$. Here, $\Th$ denotes a fine-scale mesh obtained by uniform refinement of $\TH$, where the number of refinements should be chosen such that the resulting mesh resolves all oscillations of $A$ and $f$. For solving the patch problems \cref{eq:V_H_ell_z,eq:evp}, one considers patch-local subspaces of $\V_h$.  

The SL-GFEM is straightforward to implement, as only very few technical details need to be addressed. The local spaces $\V_{H,z,h}^\ell$ (discrete counterparts of $\Vpumzi$) can be computed in parallel. Their computation only requires the local stiffness and mass matrices on the respective patches.

 In contrast to, e.g., MS-GFEM methods~\cite{BaL11,EFENDIEV2013116,Ma22}, the SL-GFEM solves local eigenvalue problems which are posed in the space spanned by a small number of deterministic snapshots. This results in a lower dimension of the eigenvalue problems~\cref{eq:evp} and, hence, makes them easier to solve numerically. After a multiplication with the respective hat-functions (partition of unity functions), we store the eigenfunctions corresponding to the~$n$ largest eigenvalues of the eigenvalue problems~\cref{eq:evp}. These functions are then used as ansatz functions for computing the global approximation~\cref{eq:PUMSLOD}. Compared to the SLOD, by construction, no stability issues in the choice of basis can occur for the SL-GFEM and thus, no special treatment of boundary patches is required, cf. \cite[App.~B]{HaPe21b}.

For the implementation, we use \texttt{gridlod} \cite{gridlod}, which is a \texttt{Python}-library initially designed for the implementation of LOD-related methods.
Although we do not require particular LOD-functionality from \texttt{gridlod}, it is convenient to use its flexible data structures for patches and its local discretization tools. Similarly, as in \cite{KR21}, our implementation can solve all local patch problems in parallel on an HPC cluster. As a comparison, we also implemented the SLOD from \cite{HaPe21b} in \texttt{gridlod}. All experiments are fully reproducible, and the corresponding source code can be found in \cite{code}\footnote{A respective \texttt{GitHub}-repository can be found in \url{https://github.com/TiKeil/SL-GFEM}}. 

\subsection{Numerical experiments}
We consider the domain $\Omega = (0,1)^2$ equipped with coarse Cartesian meshes $\TH$ and a fine Cartesian mesh $\Th$ obtained by uniform refinement of $\TH$. Note that, for ease of presentation, $H$ and $h$ henceforth denote the elements side lengths instead of their diameters. For all numerical experiments, we use $h=2^{-10}$, which results in about one million degrees of freedom for the fine mesh. Note that our implementation also works for higher spatial dimensions~$d$.
For the numerical experiments, we consider two scalar diffusion coefficients $A$ (realization of random field with short correlation length and high contrast channeled coefficient) and two source terms $f$ (constant and non-polynomial). The precise definitions can be found in the respective experiments. Each configuration serves its own purpose for numerically investigating the SL-GFEM. For measuring the approximation quality, we use the relative energy error, i.e., 
\begin{equation*}
	e_{a}^{\text{rel}}(\tilde{u}) \coloneqq \frac{\vnormof{u_h - \tilde{u}}}{\vnormof{u_h}},
\end{equation*}
where $u_h \in \V_h$ denotes the first order finite element approximation of~\cref{eq:sol} which we use as reference solution. Further, $\tilde{u}$ is a placeholder for the SL-GFEM approximation~\cref{eq:PUMSLOD} or the SLOD approximation~\cref{eq:SLOD}.

\subsubsection{Super-exponential localization}
\label{subs:exp1}

First, we investigate the localization properties of the SL-GFEM given several choices of the local approximation space size~$n$. For the choice $f \equiv 1$, the optimal order term in \cref{t:aposteriori,t:apriori} disappears and only the localization error and the approximation error $\dn$ are present. As coefficient~$A$, we consider a realization of the random field taking piecewise constant values on $\T_{2^{-8}}$, which are independent and identically distributed in the interval~$[1,100]$. This results in a maximum contrast of~$100$. We consider the fixed coarse mesh $\T_{2^{-5}}$ and the polynomial degree $p=0$.
\begin{figure}[h]
\begin{tikzpicture}

\definecolor{darkgray176}{RGB}{176,176,176}
\definecolor{lightgray204}{RGB}{204,204,204}

\begin{axis}[
width=0.52\textwidth,
legend cell align={left},
legend style={fill opacity=0.8, draw opacity=1, text opacity=1, draw=lightgray204, xshift=8pt, font=\small,
anchor=north west},
log basis y={10},
tick align=outside,
tick pos=left,
x grid style={darkgray176},
xmajorgrids,
xmin=0.85, xmax=4.15,
xlabel={\(\displaystyle \ell\)},
ylabel={\(\displaystyle e_{a}^{\text{rel}}\)},
y label style={at={(axis description cs:-0.05,0.5)}},
xtick={1,2,3,4},
xtick style={color=black},
y grid style={darkgray176},
ymajorgrids,
ymin=5.55070539481721e-08, ymax=0.398461802702745,
ymode=log,
ytick style={color=black}
]
\addplot [thick, \PUMnOne, mark=o, mark size=3, mark options={solid,fill opacity=0}]
table {%
1 0.141170783989801
2 0.0884911338595129
3 0.0679957546962469
4 0.0574795423135526
};
\addlegendentry{SL-GFEM, $n=10$, $p=0$}
\addplot [thick, \PUMnTwo, mark=o, mark size=4, mark options={solid,fill opacity=0}]
table {%
1 0.0129215124423327
2 0.00554939064111813
3 0.00403067668244442
4 0.0033632882201987
};
\addlegendentry{SL-GFEM, $n=15$, $p=0$}
\addplot [thick, \PUMnThree, mark=o, mark size=5, mark options={solid,fill opacity=0}]
table {%
1 0.0105943128788923
2 0.00173846145245774
3 0.00113554038865404
4 0.000894505793786141
};
\addlegendentry{SL-GFEM, $n=20$, $p=0$}
\addplot [thick, \PUMnFour, mark=o, mark size=6, mark options={solid,fill opacity=0}]
table {%
1 0.0105943128788923
2 3.44559621876808e-05
3 1.79481886524478e-05
4 1.33062292174343e-05
};
\addlegendentry{SL-GFEM, $n=30$, $p=0$}
\addplot [thick, \PUMnFive, mark=o, mark size=7, mark options={solid,fill opacity=0}]
table {%
1 0.0105943128788923
2 3.00024634821805e-06
3 7.7694539432222e-07
4 5.49390504821657e-07
};
\addlegendentry{SL-GFEM, $n=40$, $p=0$}
\addplot [thick, \SLODColor, mark=*, mark size=3, mark options={solid,fill=black!50, opacity=0.8}]
table {%
1 0.19442342120425
2 0.0032830778504905
3 4.01529827930707e-05
4 1.13759137875019e-07
};
\addlegendentry{SLOD}
\end{axis}

\end{tikzpicture}
	\caption{\small Localization errors of the SL-GFEM for multiple choices of $n$ and of the SLOD for a fixed coarse mesh.}
	\label{fig:exp1}
\end{figure}
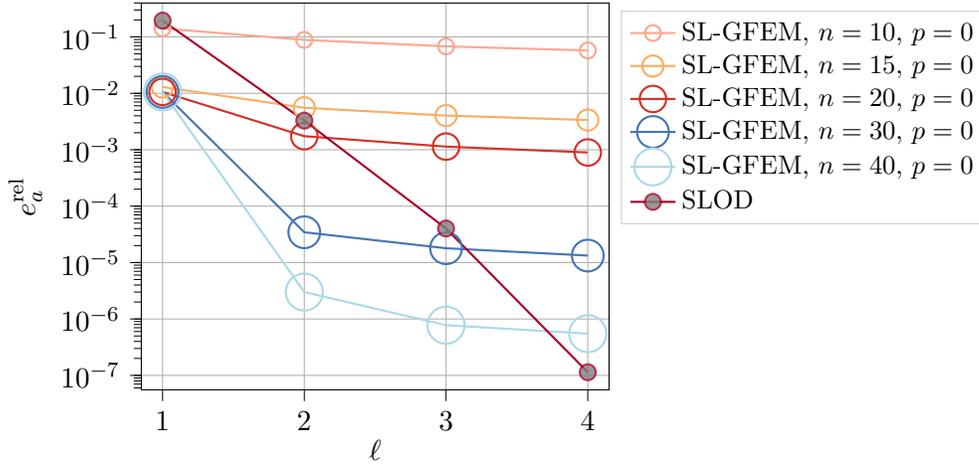
For several sizes of local approximation spaces~$n$, \cref{fig:exp1} depicts the relative energy errors of the SL-GFEM and the SLOD as a function of the oversampling parameter~$\ell$. Clearly, the parameter~$n$ strongly impacts the approximation error of the SL-GFEM.
One observes a large difference in the approximation quality, for instance, for $n=10$ and $n=15$. Conversely, choosing $n=20$ does not yield a significantly better approximation than $n=15$. This effect is related to the large jumps between the plateaus in \cref{fig:sv} and is also visible  \cref{thm:bound_delta_SLOD,thm:bound_delta_LOD}.
In conclusion, for a sufficiently large~$n$, the SL-GFEM shows a rapid decay of the localization error confirming \cref{t:aposteriori,t:apriori} numerically. For the SLOD, the super-localization property \cite[Sec.~7]{HaPe21b} is visible in \cref{fig:exp1}.

\subsubsection{Optimal convergence}
\label{subs:exp2}

\begin{figure}[h]
\begin{tikzpicture}

\definecolor{darkgray176}{RGB}{176,176,176}
\definecolor{gray127}{RGB}{127,127,127}
\definecolor{lightgray204}{RGB}{204,204,204}

\begin{axis}[
width=0.52\textwidth,
legend cell align={left},
legend style={
  fill opacity=0.8,
  draw opacity=1,
  text opacity=1,
  anchor=north west,
  draw=lightgray204,
  xshift=8pt,
  font=\small
},
log basis x={2},
log basis y={10},
tick align=outside,
tick pos=left,
x grid style={darkgray176},
xlabel={\(\displaystyle 1/H\)},
ylabel={\(\displaystyle e_{a}^{\text{rel}}\)},
y label style={at={(axis description cs:-0.05,0.5)}},
xmajorgrids,
xmin=4, xmax=34.2967508011614,
xmode=log,
xtick style={color=black},
y grid style={darkgray176},
ymajorgrids,
ymin=0.000310962603574484,  ymax=0.234314033471307,
ymode=log,
ytick style={color=black}
]
\addplot [thick, \PUMnTwo, mark=\lOneMarker, mark size=3, mark options={solid,fill opacity=0}]
table {%
4 0.0484265957365435
8 0.0214281137115357
16 0.0483333961836329
32 0.099582755032122
};
\addlegendentry{SL-GFEM, $n=10$, $p=0$, $\ell=1$}
\addplot [thick, \PUMnFour, mark=\lOneMarker, mark size=3, mark options={solid,fill opacity=0}]
table {%
4 0.0478498676414315
8 0.00895589282479858
16 0.00308194911125231
32 0.00753679075390112
};
\addlegendentry{SL-GFEM, $n=50$, $p=0$, $\ell=1$}
\addplot [thick, \PUMnTwo, mark=\lTwoMarker, mark size=5, mark options={solid,fill opacity=0}]
table {%
4 0.0407460901450893
8 0.0156175520446287
16 0.0306717654239651
32 0.0625560474233916
};
\addlegendentry{SL-GFEM, $n=10$, $p=0$, $\ell=2$}
\addplot [thick, \PUMnFour, mark=\lTwoMarker, mark size=5, mark options={solid,fill opacity=0}]
table {%
4 0.0385800971497319
8 0.00782713465253406
16 0.00179267786691774
32 0.000412421006598568
};
\addlegendentry{SL-GFEM, $n=50$, $p=0$, $\ell=2$}
\addplot [thick, \PUMnTwo, mark=\lThreeMarker, mark size=7, mark options={solid,fill opacity=0}]
table {%
4 0.0407405280735957
8 0.0142860738459989
16 0.0242682128227511
32 0.0482557156128144
};
\addlegendentry{SL-GFEM, $n=10$, $p=0$, $\ell=3$}
\addplot [thick, \PUMnFour, mark=\lThreeMarker, mark size=7, mark options={solid,fill opacity=0}]
table {%
4 0.0385662049299648
8 0.0077894231168943
16 0.00178766043154371
32 0.000410123877150939
};
\addlegendentry{SL-GFEM, $n=50$, $p=0$, $\ell=3$}
\addplot [thick, \SLODColor, mark=\lOneMarker, mark size=3, mark options={solid,fill=black!50, opacity=0.8}]
table {%
4 0.190539073433016
8 0.0558987927380846
16 0.04105460508513
32 0.0886273456256225
};
\addlegendentry{SLOD, $\ell=1$}
\addplot [thick, \SLODColor, mark=\lTwoMarker, mark size=5, mark options={solid,fill=black!50, opacity=0.8}]
table {%
8 0.0495462652902368
16 0.0119377529993168
32 0.00326564910090182
};
\addlegendentry{SLOD, $\ell=2$}
\addplot [thick, \SLODColor, mark=\lThreeMarker, mark size=7, mark options={solid,fill=black!50, opacity=0.8}]
table {%
8 0.0495458368158948
16 0.0119354650519493
32 0.00304578348800421
};
\addlegendentry{SLOD, $\ell=3$}
\addplot [semithick, color=black, dashed]
table {
	4 1e-2
	32 0.00015625
};
\addlegendentry{slope 2}
\end{axis}

\end{tikzpicture}
	\caption{\small Convergence plot of the SL-GFEM and SLOD for multiple choices of $n$ and $\ell$.}
	\label{fig:exp2}
\end{figure}
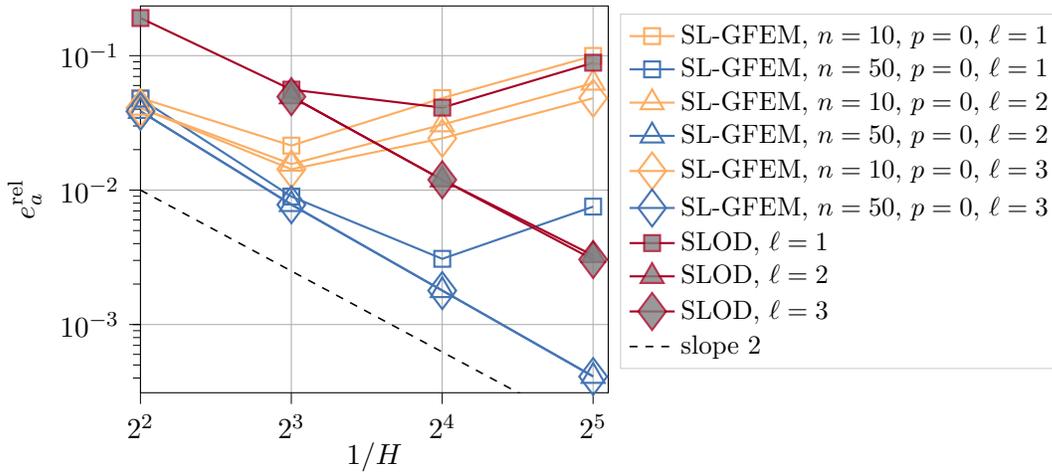

For investigating the convergence with respect to the coarse mesh size $H$, we use the same coefficient as in \cref{subs:exp1} but consider the non-polynomial source term
\begin{equation*}
	f_1(x_1,x_2) \coloneqq (x_1 + \cos(3 \pi x_1)) \cdot x_2^3\in H^1(\Omega).
\end{equation*}
\cref{fig:exp2} depicts the errors of the SL-GFEM and SLOD for multiple choices of~$n$ and~$\ell$ as a function of~$H$. As a reference, a line with slope two indicates the expected order of convergence.
For $n$ and $\ell$ sufficiently large, one observes that the SL-GFEM converges with an order of two which numerically confirms \cref{t:aposteriori,t:apriori}. Notably, the errors of the SL-GFEM are smaller by nearly one order of magnitude than the errors of the SLOD. This effect only appears for non-trivial coefficients $A$, i.e., the effect is most probably related to the contrast of the coefficient. The contrast dependence is investigated more closely in the following subsection.

\subsubsection{High-contrast channeled coefficient}
\label{subs:exp3}

\begin{figure}
	\includegraphics[width=0.48\textwidth]{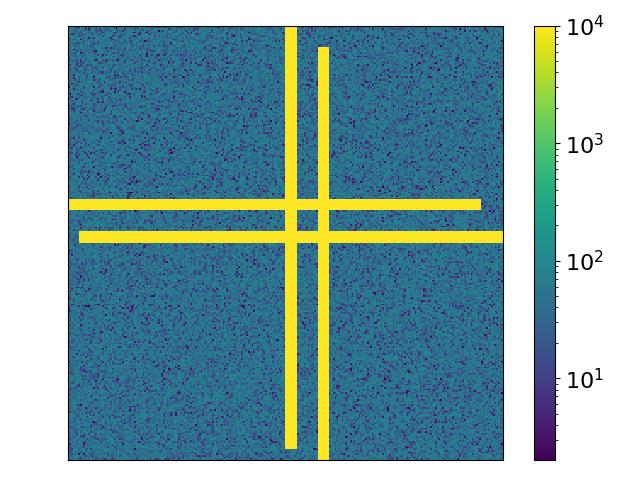}
	\includegraphics[width=0.48\textwidth]{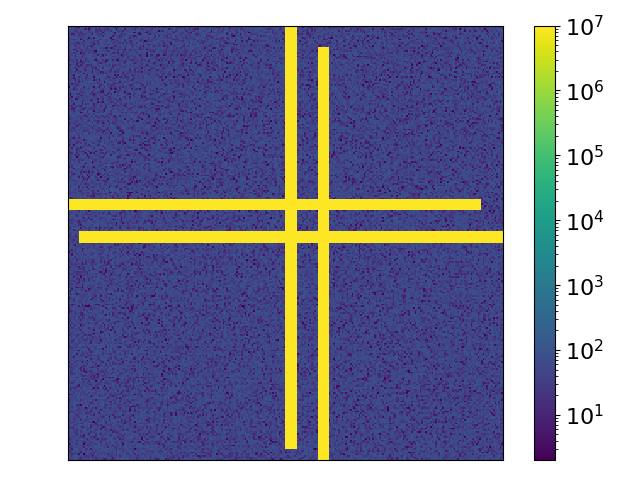}
	\caption{\small Coefficient $A_{\kappa}$ for $\kappa=10^4, 10^7$ (left and right).}
	\label{fig:exp3_coef}
\end{figure}

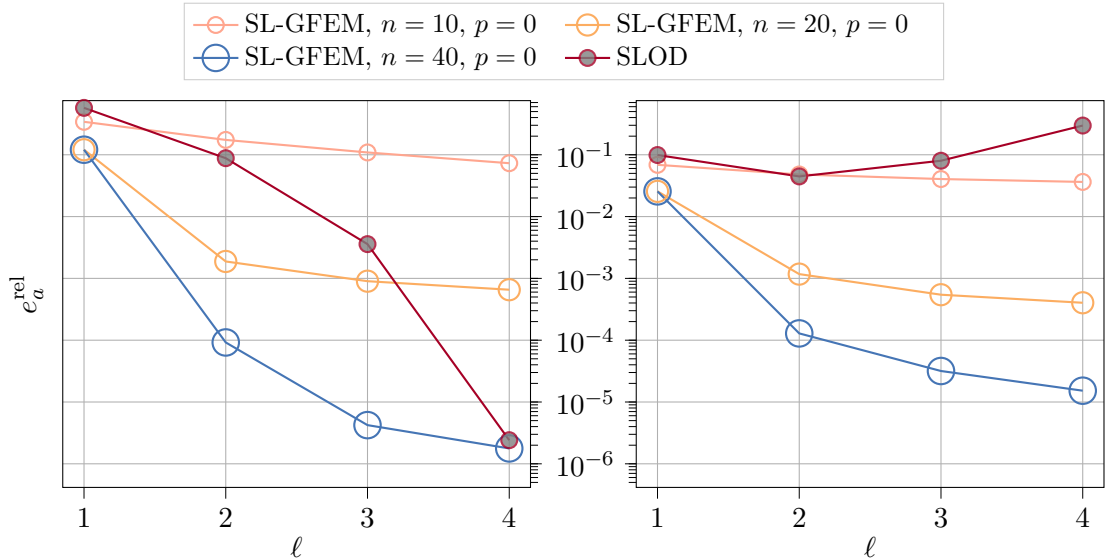
\begin{figure}[h]
\begin{tikzpicture}

\definecolor{darkgray176}{RGB}{176,176,176}
\definecolor{lightgray204}{RGB}{204,204,204}

\begin{axis}[
name=left,
width=0.52\textwidth,
legend cell align={left},
legend style={
  fill opacity=0.8,
  draw opacity=1,
  text opacity=1,
  anchor=north east,
  yshift=40pt,
  xshift=160pt,
  draw=lightgray204,
  font=\small,
  legend columns=2,
  font=\small
},
log basis y={10},
tick align=outside,
tick pos=left,
x grid style={darkgray176},
xmajorgrids,
xmin=0.85, xmax=4.15,
xlabel={\(\displaystyle \ell\)},
ylabel={\(\displaystyle e_{a}^{\text{rel}}\)},
y label style={at={(axis description cs:0.12,0.5)}},
xtick style={color=black},
xtick={1,2,3,4},
y grid style={darkgray176},
ymajorgrids,
ymin=4.15067078769619e-07, ymax=0.751660359262083,
ymode=log,
ytick style={color=black},
yticklabels={,,},
ytick pos=right,
yticklabel pos=left,
]
\addplot [thick, \PUMnOne, mark=o, mark size=3, mark options={solid,fill opacity=0}]
table {%
1 0.341416661770608
2 0.173939835390188
3 0.10901339695524
4 0.07302412554349
};
\addlegendentry{SL-GFEM, $n=10$, $p=0$ \hphantom{0}}
\addplot [thick, \PUMnTwo, mark=o, mark size=4, mark options={solid,fill opacity=0}]
table {%
1 0.120713599052476
2 0.0018778646799998
3 0.000900844856825033
4 0.00065551415092622
};
\addlegendentry{SL-GFEM, $n=20$, $p=0$ \hphantom{0}}
\addplot [thick, \PUMnFour, mark=o, mark size=5, mark options={solid,fill opacity=0}]
table {%
1 0.120713599052476
2 9.1945696373516e-05
3 4.22556592988036e-06
4 1.77410721083959e-06
};
\addlegendentry{SL-GFEM, $n=40$, $p=0$ \hphantom{0}}
\addplot [thick, \SLODColor, mark=*, mark size=3, mark options={solid,fill=black!50,opacity=0.8}]
table {%
1 0.573933387488881
2 0.0880049512946704
3 0.00358150580831056
4 2.41111469762332e-06
};
\addlegendentry{SLOD}
\end{axis}

\begin{axis}[
name=right,
at=(left.south east),
anchor=south west,
xshift=40pt,
width=0.52\textwidth,
legend cell align={left},
log basis y={10},
tick align=outside,
tick pos=left,
xlabel={\(\displaystyle \ell\)},
x grid style={darkgray176},
xtick={1,2,3,4},
xmajorgrids,
xmin=0.85, xmax=4.15,
xtick style={color=black},
y grid style={darkgray176},
ymajorgrids,
ymin=4.15067078769619e-07, ymax=0.751660359262083,
ymode=log,
ytick style={color=black},
]
\addplot [thick, \PUMnOne, mark=o, mark size=3, mark options={solid,fill opacity=0}]
table {%
1 0.0687234230035897
2 0.0480209526577795
3 0.0404890705697322
4 0.0363570031019531
};
\addplot [thick, \PUMnTwo, mark=o, mark size=4, mark options={solid,fill opacity=0}]
table {%
1 0.0255693276679326
2 0.00117268235855824
3 0.000544358824742104
4 0.000403271751422638
};
\addplot [thick, \PUMnFour, mark=o, mark size=5, mark options={solid,fill opacity=0}]
table {%
1 0.0255693276679326
2 0.000128783160458677
3 3.17361119559097e-05
4 1.52372658506678e-05
};
\addplot [thick, \SLODColor, mark=*, mark size=3, mark options={solid,fill=black!50,opacity=0.8}]
table {%
1 0.0989749371197156
2 0.0446295389155419
3 0.0801053765718436
4 0.297055395894572
};
\end{axis}

\end{tikzpicture}
	\caption{\small Localization errors of the SL-GFEM and the SLOD for the high-contrast channeled coefficient $A_\kappa$ for $\kappa = 10^4, 10^7$ (left and right).}
	\label{fig:exp3}
\end{figure}

One of the major challenges for multi-scale methods is their sensitivity to high-contrast channeled coefficients.  
In this numerical experiment, we consider the coefficient $A_{\kappa}$ constructed by adding four channels of conductivity~$\kappa$ in the coefficient from \cref{subs:exp1,subs:exp2}. Some of the channels touch the boundary, while others stop before. The number $\kappa$ is the maximum contrast of~$A_\kappa$. \Cref{fig:exp3_coef} illustrates the coefficient for $\kappa = 10^5, 10^7$.
For this numerical experiment, we choose the same setup as in \cref{subs:exp1}. \Cref{fig:exp3} depicts the localization errors for the above choices of $\kappa$.
The SL-GFEM appears to be largely unaffected by large values of~$\kappa$. For the SLOD, in contrast, the best practical realization known until now~\cite{HaPe21b} yields a basis with deteriorating stability as $\kappa$ is increased (growing constants in~\cref{a:riesz,a:rieszloc}). This explains the worse performance of the SLOD for $\kappa = 10^7$ compared to $\kappa = 10^4$. 
Notably, when compared to \cref{subs:exp1}, the SL-GFEM does not need more local functions~$n$ to attain a good approximation quality, which suggests that the choice of~$n$ is not affected by the contrast.

\subsubsection{Higher-order polynomials}
\label{subs:exp4}

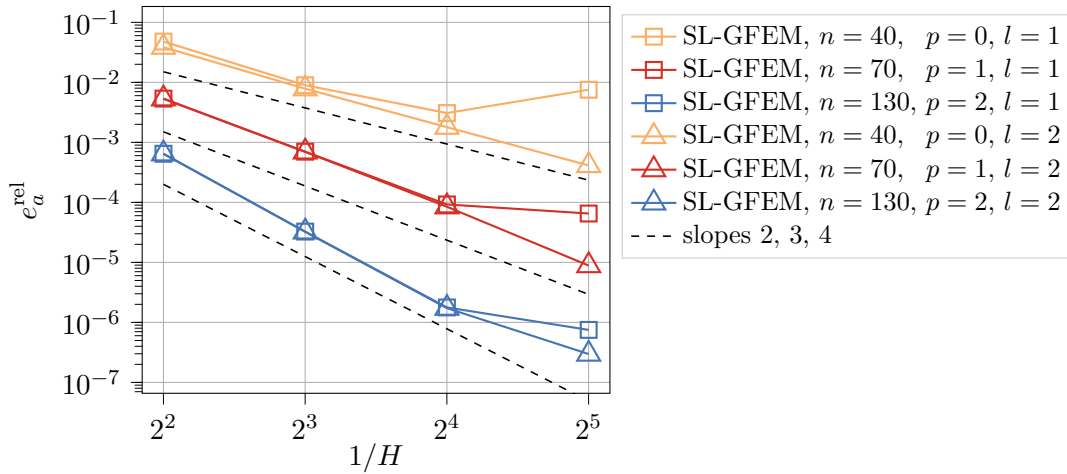
\begin{figure}[h]
\begin{tikzpicture}

\definecolor{darkgray176}{RGB}{176,176,176}
\definecolor{lightgray204}{RGB}{204,204,204}

\begin{axis}[
width=0.52\textwidth,
legend cell align={left},
legend style={fill opacity=0.8, draw opacity=1, text opacity=1, draw=lightgray204,
xshift=8pt, anchor=north west, font=\small},
log basis x={2},
log basis y={10},
tick align=outside,
tick pos=left,
x grid style={darkgray176},
xlabel={\(\displaystyle 1/H\)},
ylabel={\(\displaystyle e_{a}^{\text{rel}}\)},
y label style={at={(axis description cs:-0.05,0.5)}},
xmajorgrids,
xmin=3.60500185044332, xmax=35.506223106171,
xmode=log,
xtick style={color=black},
y grid style={darkgray176},
ymajorgrids,
ymin=6.58794904368016e-08, ymax=0.183720237778836,
ymode=log,
ytick style={color=black}
]
\addplot [thick, \PUMnTwo, mark=\lOneMarker, mark size=3, mark options={solid,fill opacity=0}]
table {%
4 0.0478498676414315
8 0.00895589282479858
16 0.00308194911125231
32 0.00753679075390112
};
\addlegendentry{SL-GFEM, $n=40$,\hphantom{0} $p=0$, $l=1$}
\addplot [thick, \PUMnThree, mark=\lOneMarker, mark size=3, mark options={solid,fill opacity=0}]
table {%
4 0.00536205097371743
8 0.000699393352387509
16 9.2391908819023e-05
32 6.52161826706634e-05
};
\addlegendentry{SL-GFEM, $n=70$,\hphantom{0} $p=1$, $l=1$}
\addplot [thick, \PUMnFour, mark=\lOneMarker, mark size=3, mark options={solid,fill opacity=0}]
table {%
4 0.000645466310783716
8 3.26046474792041e-05
16 1.77610687672464e-06
32 7.49632117467066e-07
};
\addlegendentry{SL-GFEM, $n=130$, $p=2$, $l=1$}
\addplot [thick, \PUMnTwo, mark=\lTwoMarker, mark size=5, mark options={solid,fill opacity=0}]
table {%
4 0.0385800971497319
8 0.00782713465253406
16 0.00179267786691774
32 0.000412421006598568
};
\addlegendentry{SL-GFEM, $n=40$,\hphantom{0} $p=0$, $l=2$}
\addplot [thick, \PUMnThree, mark=\lTwoMarker, mark size=5, mark options={solid,fill opacity=0}]
table {%
4 0.00532514550154908
8 0.000693839005442121
16 8.45293775765913e-05
32 8.84717371274897e-06
};
\addlegendentry{SL-GFEM, $n=70$,\hphantom{0} $p=1$, $l=2$}
\addplot [thick, \PUMnFour, mark=\lTwoMarker, mark size=5, mark options={solid,fill opacity=0}]
table {%
4 0.000645658886983883
8 3.22445914819123e-05
16 1.72647362023178e-06
32 2.96748098275482e-07
};
\addlegendentry{SL-GFEM, $n=130$, $p=2$, $l=2$}
\addplot [semithick, color=black, dashed]
table {
	4 1.5e-2
	32 0.00023437499999999988
};
\addlegendentry{slopes 2, 3, 4}
\addplot [semithick, color=black, dashed, forget plot]
table {
	4 1.5e-3
	32 2.9296875e-06
};
\addlegendentry{slope 3}
\addplot [semithick, color=black, dashed, forget plot]
table {
	4 2e-4
	32 4.882812500000008e-08
};
\addlegendentry{slope 4}
\end{axis}

\end{tikzpicture}
	\caption{\small Convergence plot of the higher-order SL-GFEM for multiple choices of $n$, $\ell$, and $p$.}
	\label{fig:exp4}
\end{figure}

One key benefit of the proposed method is its flexibility with regard to the choice of polynomial degree, i.e., the construction of higher-order methods is straightforward. While the previous numerical experiments have investigated the performance of the SL-GFEM for $p = 0$, this experiment also considers higher polynomial degrees. Using the setup from \cref{subs:exp2}, \cref{fig:exp4} depicts the errors of the SL-GFEM for $p = 0,1,2$ as a function of $H$ together with lines indicating the respective expected orders of convergence.
For $n$ and $\ell$ sufficiently large, it can be observed that the method of degree $p$ converges with an order of $p+2$ (recall that $f$ is sufficiently smooth). This numerically confirms \cref{t:aposteriori,t:apriori}. 
Note that the choice of~$n$ needs to be adapted to~$p$, which is related to the larger plateaus in \cref{fig:sv}. We observed that $n$ needs to be increased linearly as $p$ is increased. It is left to future research to find a (possibly adaptive) choice of $n$ such that pessimistic choices (that may result in unnecessarily many basis functions) can be avoided.

\bibliographystyle{alpha}
\bibliography{bib}
\end{document}